\documentclass[12pt,leqno]{article}
\usepackage{amssymb}
\usepackage{amsmath}
\usepackage{natbib}
\begin{document}
\bibliographystyle{plainnat}
 \oddsidemargin= 20mm \hoffset=-4pc \topmargin= 35mm \textwidth=175mm
\textheight=270mm \pagestyle{plain}
%\newcounter{r}
%\newcounter{st}
%\newcounter{et}[st]
%\newcommand{\en}{\thest.\refstepcounter{et}}
%\numberwithin{equation}{section}
\renewcommand{\theequation}{\thesection.\arabic{equation}}
\newcommand{\sign}{sign}
\newcommand {\Log}{Log}
\newcommand{\Res}{Res}
\title{ Title: An expansion of zeta(3) in continued fraction with parameter.}

\author{Authors:L.A.Gutnik}

\date{}
\maketitle

Comments: 40 pages with 0 figures

Subj-class: Number Theory

MSC-class: 11-XX \vskip.20pt
 \begin{abstract} We present here continued
fraction for Zeta(3) parametrized by some family of points (F,G) on projective
line. This family of points can be obtained if  from full projective line would
be removed some no more than countable nowhere dense
  exeptional set of finite points.
 A countable nowhere dense  set, which contains the above exeptional set of
 finite points, is specified also.
 \end{abstract} \vskip.20pt
%\hspace*{7cm}
 {\hfill\sl To the thirty fifth anniversary }

\hspace*{8cm}{\hfill\sl of Ap\'ery's discovery.}
\vskip.15pt

%%%%%%%%%%%%%%%%%%%%%%%%%%%%%%%% Table of contents %%%%%%%%%%%%%%%%%%%%%%%%%%%%%%%%%%%%%%
\begin{center}{\large\bf Table of contents}\end{center}
\vskip.5pt \S 0. Foreword. \vskip.5pt\noindent \S 1. Introduction. Begin of the
proof of Theorem B. \vskip.5pt\noindent \S 2. Transformation of the system
considered in \S 1. \vskip.5pt\noindent \S 3. Calculation of the matrix
$A^{\ast\ast}_{1,0}(z,\nu).$ \vskip.5pt\noindent \S 4. Further properties of
the functions considered in \S 1. \vskip.5pt\noindent \S 5. Auxiliary
difference equation. \vskip.5pt\noindent \S 6. Auxiliary continued fraction.
\vskip.5pt\noindent \S 7. End of the Proof of theorem B.
%%%%%%%%%%%%%%%%%%%%%%%%%%%%%%% 6.0%%%%%%%%%%%%%%%%%%%%%%%%%%%%%%%%%%%%%%
{\begin{center}\large\bf\S 0. Foreword . \end{center}}
\newpage
\pagestyle{headings}
\topmargin= -0.7mm \textheight=236mm \vskip.6mm \markright{\footnotesize\bf
L.A.Gutnik, An expansion of zeta(3) in continued fraction with parameter.}
 We say that two infinite continuous fractions
are equivalent if the set of their common convergents is infinite.
We say that two infinite continuous fractions
are essentially distinct if the set of their common convergents
 is finite or empty.

If $\delta_0=1,\delta_\nu\in{\mathbb C}\diagdown\{0\},$
$$b_\nu^{(2)}=b_\nu^{(1)}\delta_\nu,
a_\nu^{(2)}=a_\nu^{(1)}\delta_\nu\delta_{\nu-1}$$
for $\nu\in{\mathbb N},$
 then easy induction show that
$P_\nu^{(2)}=P_\nu^{(1)}\prod\limits_{\kappa=1}^\nu\delta_\kappa,
Q_\nu^{(2)}=Q_\nu^{(1)}\prod\limits_{\kappa=1}^\nu\delta_\kappa,$
and  continued fraction
\begin{equation}\label{eq:ithcf}
b_0^{(i)}+\frac{a_1^{(i)}\vert}{\vert b_1^{(i)}}+
\frac{a_2^{(i)}\vert}{\vert b_2^{(i)}}+...+
\frac{a_\nu^\vee\vert}{\vert b_\nu^{(i)}}+...,
\end{equation}
with i=1 is equivalent
 to continued fraction (\ref{eq:ithcf}) with $i=2.$
 According to the famous result of R. Ap\'ery  \citep{Apery},
\begin{equation}\label{eq:0a11}
\zeta(3)=
b_0^\vee+\frac{a_1^\vee\vert}{\vert b_1^\vee}+
\frac{a_2^\vee\vert}{\vert b_2^\vee}+...+
\frac{a_\nu^\vee\vert}{\vert b_\nu^\vee}+...
\end{equation}
with
\begin{equation}\label{eq:g2}
b_0^\vee=0,\,b_1^\vee=5,\,a_1^\vee=6,\,b_{\nu+1}^\vee=
\end{equation}
$$34\nu^3+51\nu^2+27\nu+5,\,a_{\nu+1}^\vee=-\nu^6,$$
for $\nu\in{\mathbb N}.$
We denote by $r_{A,\nu}$ the $\nu$-th  convergent of
continued fraction (\ref{eq:0a11}).

\noindent Yu.V. Nesterenko in \citep{Nest} has offered the following expansion
of $2\zeta(3)$ in continuous fraction:
\begin{equation}\label{eq:N1}
2\zeta(3)=2+\frac{1\vert}{\vert2}+\frac{2\vert}{\vert4}+\frac{1\vert}{\vert3}+
\frac{4\vert}{\vert2}...,
\end{equation}
with
\begin{equation}\label{eq:N2}
b_0=b_1=a_2=2,\,a_1=1,\,b_2=4,
\end{equation}
\begin{equation}\label{eq:N3}
b_{4k+1}=2k+2,\,a_{4k+1}=k(k+1),\,b_{4k+2}=
\end{equation}
$$2k+4,\,a_{4k+2}=(k+1)(k+2)$$
for $k\in{\mathbb N},$
\begin{equation}\label{eq:N4}
b_{4k+3}=2k+3,\,a_{4k+3}=(k+1)^2,\,b_{4k+4}=
\end{equation}
$$2k+2,\,a_{4k+4}=(k+2)^2\,\,\text{for}\,\, k\in{\mathbb N}_0.$$
\newpage We denote by $r_{N,\nu}$ the $\nu$-th convergent of continued
fraction (\ref{eq:N1}). The continued fractions (\ref{eq:N1}) and
(\ref{eq:0a11}) are equivalent because $2r_{A,\nu}=r_{N,4\nu-2}$ for all
$\nu\in{\mathbb N}.$

Elementary proof of Yu.V. Nesterenko result

 can be found in \citep{G16}.

Let me to formulate my result now. Let $u$ and $v$ are variables,
% $F$ and $G$ are constants in ${\mathbb R}$ such that either
%\begin{equation}\label{eq:FG}
%\vert F \vert+\vert G \vert>0,
%\end{equation}
$$\tau=\tau(\nu)=\nu+1,\,\sigma=\sigma(\nu)=\tau(\tau-1)=\nu(\nu+1),$$
where $\nu\in{\mathbb N}.$ Let further

\begin{equation}\label{eq:cfg2n}
c_{u,v,2}(\nu)=-\tau(\tau+1)^2(2\tau-1)\times
\end{equation}
$$(-3(2\tau-1)^2u^2-(10\tau^2-10\tau+3)uv+2(3\tau^4-6\tau^3+4\tau^2-\tau)v^2=$$
$$-3(4\sigma(\nu)+1)u^2-(10\sigma(\nu)+3)uv+
2\sigma(3\sigma(nu)+1)v^2\in{\mathbb N}[u,v].$$
\begin{equation}\label{eq:cfg1n}
c_{u,v,1}(\nu)=-12(68\tau^6-45\tau^4+12\tau^2-1)u^2-\end{equation}
$$8(157\tau^6-106\tau^4+30\tau^2-3)uv+$$
$$4(102\tau^8-170\tau^6+89\tau^4-24\tau^2+3)v^2\in{\mathbb N}[u,v],$$
\begin{equation}\label{eq:cfg0n}
c_{u,v,0}(\nu)=\tau(\tau-1)^2(2\tau+1)\times\end{equation}
$$(3(4\tau^2+4\tau+1)u^2+(10\tau^2+10\tau+3)uv
-2(3\tau^4+6\tau^3+4\tau^2+\tau)v^2)=$$
$$\tau(\tau-1)^2(2\tau+1)\times$$
$$(3(4\sigma(\nu+1)+1)u^2+(10\sigma(\nu+1)+3)uv
-2\sigma(\nu+1)(3\sigma(\nu+1)+1)v^2).$$
\begin{equation}\label{eq:bfgnp1}
b_{u,v}(\nu+1)=-c_{u,v,1}(\nu)\in{\mathbb Q}[u,v]\,\,\text{for}\,\,
\nu\in{\mathbb N},
\end{equation}
\begin{equation}\label{eq:afgnp1}
a_{u,v}(\nu+1)=-c_{u,v,0}(\nu)c_{u,v),2}(\nu-1)\,\,\text{for}\,\,
\nu\ge2,\,\nu\in{\mathbb N},
\end{equation}
\begin{equation}\label{eq:avg2}
a_{u,v}(2)=-c_{u,v,0}(1),
\end{equation}
\begin{equation}\label{eq:bfg0}
P_{u,v}(0)=b_{u,v}(0)=4(3u+2v),\,Q_{u,v}(0)=1,
\end{equation}
\begin{equation}\label{eq:bfg1}
Q_{u,v}(1)=b_{u,v}(1)=(34u+52v)/(u+v)
\end{equation}
\begin{equation}\label{eq:Pfg1}
P_{u,v}(1)=(327u+500v)
\end{equation}
\begin{equation}\label{eq:afg1}
a_{u,v}(1)=P_{u,v}(1)=-b_{u,v}(0)b_{u,v}(1).
\end{equation}
Calculations described in  \citep{G18} lead to the following continued fraction
over the field ${\mathbb Q}(u,v):$
\begin{equation}\label{eq:G0}
b_{u,v}(0)+\frac{a_{u,v}(1)\vert}{\vert b_{u,v}(1)}+
\frac{a_{u,v}(2)\vert}{\vert b_{u,v}(2)}+
\frac{a_{u,v}(3)\vert}{\vert b_{u,v}(3)}+... .
\end{equation}
We denote by $r_{u,v}(\nu)$
  the $\nu$-th convergent of continuous fraction
(\ref{eq:G0}). We denote by $P_{u,v}(\nu)$ and $Q_{u,v}(\nu)$
 the respectively nominator and denominator of  $r_{u,v}(\nu).$
Let further $\Delta(x)=18x^2(2x+1)+(7x+3)^2,$
$$\rho_k(x)=\frac{5x+3+(-1)^k\sqrt{\Delta(x)}}{x(3x+2)},$$
where $x=2\nu(\nu+1),\,\nu\in{\mathbb N},\,k=1,2,$ and let
$${\mathfrak A}=\{\rho_k(x)\colon
x=2\nu(\nu+1),\,\nu\in{\mathbb N},\,k=1,2\}.$$
$$
 {\mathfrak B}=\bigg\{
-\beta^{\ast(2)}_2(1;\nu)/\beta^{\ast(1)}_2(1;\nu) \colon \nu\in{\mathbb
N}_0\bigg\}, $$
 where
\begin{equation}\label{eq:barzn}
\beta^{\ast
(r)}_2(z;\nu)=\sum\limits_{k=0}^{\nu+1}\left(\binom{\nu+1}k\binom{\nu+k}k\right)^2t^rz^k.
\end{equation}

 Then we have
 
\bf Theorem B. \it Let
 $F+G\ne0,\,(F+G)G\ge0,\,G/F\not\in{\mathfrak B}\cup{\mathfrak A},$
 if $F\ne0.$
 Then the specialization of (\ref{eq:G0}) for $u=F,\,v=G$ is well
defined (i.e all convergents $r_{u,v}(\nu)$ are well defined
 for $u=F,\,v=G$) and the following equality holds
\begin{equation}\label{eq:G1}
8(F+G)\zeta(3)=b_{F,G}(0)+\frac{a_{F,G}(1)\vert}{\vert b_{F,G}(1)}+
\frac{a_{F,G}(2)\vert}{\vert b_{F,G}(2)}+ \frac{a_{F,G}(3)\vert}{\vert
b_{F,G}(3)}+... .
\end{equation}
moreover $P_{u,v}(\nu)$ and $Q_{u,v}(\nu)$ are homogeneous polynomials in
${\mathbb Z}[u,v],$ and
\begin{equation}\label{eq:degP}
\max(2\nu,1)=\deg_u(P_{u,v}(\nu))=\deg_v(P_{u,v}(\nu))=\deg(P_{u,v}(\nu)),
\end{equation}
\begin{equation}\label{eq:degQ}
\max(2\nu-1,0)=\deg_u(Q_{u,v}(\nu))=\deg_v(Q_{u,v}(\nu))=\deg(Q_{u,v}(\nu)),
\end{equation}
where $\nu\in {\mathbb N}_0.$

\bf Remark. \rm The values $\rho_k(x)$ with $k=1,2$
 are zeros of the following trinomial $a_0(x)+2a_1(x)\rho+a_2(x)\rho^2,$
 where $$a_0(x)=-6(2x+1),\,a_1(x)=-2(5x+3),\,a_2(x)=x(3x+2).$$
Since $-a_0(x)/a_2(x)=3/x+3/(3x+2), -a_1(x)/a_2(x)=3/x+1/(3x+2),$ decrease
together with increasing of $x>0$ it follows that
$$\Delta(x)/(a_2(x))^2=(-a_1(x)/a_2(x))^2-
(-a_0(x)/a_2(x)),$$ and $\rho_2(x)$ decrease with increasing of $x.$
 Moreover, $0<-\rho_1(x)<\rho_2(x)$
for any $x>0$ and $\lim\limits_{x\to\infty} r_2(x)=0.$ Consequently, for given
$F>0$ and $G>0$ the condition of the Theorem B
 must be checked for finite family of $\nu;$ for example, if $G/F>r(4),$
 then condition of the Theorem B is fulfilled. We note that $r(4)<0,36.$

I prove Theorem B in sections 1 -- 7.
 Initial variants of this article can be found in \citep{G17},\,
 \citep{G18}.
%Let $T(u,v)\in K[u,v],$ where $K$ is a field. If $T(u,v)$ is zero
%polynomial, we put $\deg(T(u,v)=-\infty.$ Then standard relations
%$$
%\deg(T_1T_2)=\deg(T_1)+\deg(T_2),\,
%\deg(T_1+T_2)=\le \max(\deg(T_1),\deg(T_2))
%$$
%hold.

%%%%%%%%%%%%%%%%%%%%%%%%%%%%%%%%%%%%%% \S1%%%%%%%%%%%%%%%%%%%%%%%%%%%%%%%
\refstepcounter{section}
{\begin{center}\large\bf\S 1. Introduction. Begin of the proof of Theorem B.
\end{center}}
Let
\setcounter{equation}0
\begin{equation}\label{eq:1a}
\vert z\vert>1,-3\pi/2<\arg(z)\le\pi/2,\log(z)=\ln(\vert z\vert)+i\arg(z).
\end{equation}
Clearly, $\log(-z)=\log(z)-i\pi,$ when $\Re(z)>0,$
 and $\log(z)=\log(-z)-i\pi,$ when $\Re(z)<0.$
Let
\begin{equation}\label{eq:b}
f_1^{\ast}(z, \nu):=
\sum\limits_{k=0}^{\nu+1}(z)^k\binom{\nu+1}k^2\binom{\nu+k}\nu^2,
\end{equation}
\begin{equation}\label{eq:d}
R(t,\nu)=\left(\prod\limits_{j=1}^\nu(t-j)\right){\bigg /}
\left(\prod\limits_{j=0 }^{\nu+1}(t+j)\right),
\end{equation}
\begin{equation}\label{eq:h}
f_2^{\ast}(z, \nu)=
\sum\limits_{t=1}^{+\infty}z^{-t}
(\nu+1)^2(R(\alpha,t,\nu))^2,
\end{equation}

\begin{equation}\label{eq:ab2}
f_4^\ast(z,\nu)=-\sum\limits_{t=1}^{+\infty} z^{-t}
(\nu+1)^2\left (\frac\partial{\partial t}
(R^2)\right)(t,\nu),
\end{equation}

\begin{equation}\label{eq:ac}
f_3^{\ast}(z,\nu)=(\log(z))f_2^\ast(z,\nu)+
f_4^{\ast}(z,\nu),
\end{equation}

\begin{equation}\label{eq:acx}
f_k(z,\nu)=f_k^{\ast}(z,\nu)/(\nu+1)^2
\end{equation}
where $k=1,\,2,\,3,\,4,\,\nu\in{\mathbb N}_0.$
Let
\begin{equation}\label{eq:93bb}
\tau=\tau(\nu)=\nu+1,\,\mu=\mu(\nu)=\tau^2=(\nu+1)^2,
\end{equation}
\begin{equation}\label{eq:a11a}%!
a_{1,1}^\ast(z,\nu)=\frac12(-5\mu+3\mu-6\mu^2)-z\mu(1+18\mu)+
\end{equation}
$$(3\mu+\mu^2+z(7\mu+16\mu^2))\tau,$$
\begin{equation}\label{eq:a12a}%!
a_{1,2}^\ast(z;\nu)=-8\mu-4\mu^2-z(12\mu+20\mu^2)+(2+10\mu+z(2+26\mu))\tau,
\end{equation}
\begin{equation}\label{eq:a13}%!
a_{1,3}^\ast(z;\nu)=-1-14\mu+z(-1+2\mu)+(7+8\mu+z(3-8\mu))\tau,
\end{equation}
\begin{equation}\label{eq:a14a}%!
a_{1,4}^\ast(z;\nu)=(z-1)(2+12\mu-10\tau),
\end{equation}
\begin{equation}\label{eq:a21a}%!
a_{2,1}^\ast(z;\nu)=-z\mu-z20\mu^2-z12\mu^3+\tau(7z\mu+26z\mu^2),
\end{equation}
$$z\mu(24-22\alpha+5\alpha^2+28\mu-2\alpha\mu),$$
\begin{equation}\label{eq:a22a}%!
a_{2,2}^\ast(z;\nu)=-\mu-3\mu^2+z\mu(-13-38\mu)+
\end{equation}
$$(3\mu+\mu^2+2z+33z\mu+16z\mu^2)\tau,$$
\begin{equation}\label{eq:a23a}%!
a_{2,3}^\ast(z;\nu)=-8\mu-4\mu^2-z-6z\mu+4z\mu^2+(2+10\mu+5z-2z\mu)\tau
\end{equation}
\begin{equation}\label{eq:a24a}%!
a_{2,4}^\ast(z;\nu)=(1+14\mu-7\tau-8\mu\tau)(z-1),
\end{equation}
\begin{equation}\label{eq:a31a}%!
a_{3,1}^\ast(z;\nu)=-z\mu-21z\mu^2-26z\mu^3+(7z\mu+33z\mu^2+8z\mu^3)\tau,
\end{equation}
\begin{equation}\label{eq:a32a}%!
a_{3,2}^\ast(z;\nu)=-14z\mu-58z\mu^2-12z\mu^3+(2z+40z\mu+42z\mu^2)\tau,
\end{equation}
\begin{equation}\label{eq:a33a}%1
a_{3,3}^\ast(z;\nu)=
-\mu-3\mu^2-z-17z\mu-6z\mu^2+(3\mu+\mu^2+7z+17z\mu)\tau
\end{equation}
\begin{equation}\label{eq:a34a}%!
a_{3,4}^\ast(z;\nu)=(8\mu+4\mu^2-2\tau-10\mu\tau)(z-1),
\end{equation}
\begin{equation}\label{eq:a41a}%!
a_{4,1}^\ast(z;\nu)=-z\mu-21z\mu^2-38z\mu^3+(7z\mu+35z\mu^2+18z\mu^3)\tau,
\end{equation}
\begin{equation}\label{eq:a42a}%!
a_{4,2}^\ast(z;\nu)=-15z\mu-79z\mu^2-38z\mu^3+
\end{equation}
$$(2z+47z\mu+75z\mu^2+8z\mu^3)\tau,$$
\begin{equation}\label{eq:a43a}%!
a_{4,3}^\ast(z;\nu)=-z-31z\mu-48z\mu^2-4z\mu^3+(9z+53z\mu+22z\mu^2)\tau,
\end{equation}
\begin{equation}\label{eq:a44a}%!
a_{4,4}^\ast(z;\nu)=-\mu-3\mu^2-z-9z\mu-2z\mu^2+(3\mu+\mu^2+5z+7z\mu)\tau
\end{equation}
We denote by $A^\ast(z;\nu)$ the $4\times4$-matrix with $a_{i,k}^\ast(z;\nu)$
in its $i$-th row and $k$-th column for $i=1,\,...,\,4,\,k=1,\,...,\,4.$
Clearly,
\begin{equation}\label{eq:15d}
A^\ast(z;\nu)=
A^\ast(1;\nu)+(z-1)V^\ast(\nu),
\end{equation}
where the matrix $V^\ast(\nu)$ does not depend from $z.$
Let
\begin{equation}\label{eq:15e}
X_k(z;\nu)=
\left(\begin{matrix}
f_k(z,\nu)\\
\delta f_k(z,\nu)\\
\delta^2 f_k(z,\nu)\\
\delta^3 f_k(z,\nu)
\end{matrix}\right)
,\,X_k^\ast(z;\nu)=(\nu+1)^2
X_k(z;\nu)
\end{equation}
for
$k=1,\,2,\,3,\,\vert z\vert>1,\,\nu\in{\mathbb N}_0.$
Let further
\begin{equation}\label{eq:93be0}
X_k(z;-\nu-2)=X_k(z;\nu),
\end{equation}
where $\nu\in{\mathbb N}_0.$
Let us consider the row
\begin{equation}\label{eq:113bd1}
R(\nu)=(r_1(\nu),\,r_2(\nu),\,r_3(\nu),\,r_4(\nu)),
\end{equation}
where
\begin{equation}\label{eq:93bd}
r_1(\nu)=\mu(\nu)^2,\,r_2(\nu)=0,\,r_3(\nu)=-2\mu(\nu),\,r_4(\nu)=0.
\end{equation}
We have the following equalities:
$$A^\ast(z;\nu)=A_{1,0}^\ast(z;\nu),\,X_k(z;\nu)=X_{1,0,k}(z;\nu),$$
$$R(\nu)=R_{1,0}(\nu),$$
 where $A_{\alpha,0}^\ast(z;\nu),\,X_{\alpha,0,k}(z;\nu)$

and $R_{\alpha,0}(\nu)$ are studied in \citep{G5} $-$ \citep{G15}. We take
$\alpha=1$ in (105), \citep{G13}, in\ \ (1), \citep{G15}, in \S 10.1, \citep
{G14},\ \ \S 11.3,\citep{G15}. Then we have the following Theorem:

\bf Theorem 1. \it The column $X_k(z;\nu)$ satisfies
to the equation
\begin{equation}\label{eq:93ca}
\nu^5X_k(z;\nu-1)=
A^\ast(z;\nu)X_k(z;\nu),
\end{equation}
for
$\nu\in M_1^\ast=(-\infty,-2]\cup[1,+\infty))\cap{\mathbb Z},
\,k=1,\,2,\,3,\,\vert z\vert>1;$ moreover,
 the matrix $A^\ast(z;\nu)$ has the following property:
\begin{equation}\label{eq:93cc}
-\nu^5(\nu+1)^5E_{4}=
A^\ast(z;-\nu-1)A^\ast(z;\nu),
\end{equation}
where $E_4$ is the $4\times4$ unit matrix, $z\in{\mathbb C},\,\nu\in{\mathbb
C}.$ \rm The Lemma 11.3.1 in \citep{G15} have the following
 formulation for $\alpha=1$ :

\bf Theorem 2. \it The row $R(\nu)$
has the following property:
\begin{equation}\label{eq:15f}
R(\nu-1)A^\ast(1;\nu)=\nu^5R(\nu),\,\text{where}\,\,\nu\in{\mathbb C}.
\end{equation}\rm
%%%%%%%%%%%%%%%%%%%%%%%%%%%%%\S2%%%%%%%%%%%%%%%%%%%%%%%%%%%%%%%%%%%%%%%%%%%
\refstepcounter{section}
%{\begin{center}\large\bf\S 2. Transformations of the system considered in
%the Introduction.\end{center}}
%in the case $\alpha=1.$
{\begin{center}\large\bf\S 2. Transformation of the system considered
 in \S 1.\end{center}}
In view of (\ref{eq:93bb}), (\ref{eq:93bd})
\setcounter{equation}0
\begin{equation}\label{eq:15.1.a}
r_1(\nu)=\mu(\nu)^2=(\nu+1)^4=\tau^4,\,r_2(\nu)=0,\,r_3(\nu)=
\end{equation}
$$-2\mu(\nu)=-2(\nu+1)^2=-2\tau^2,\,r_4(\nu)=0.$$
Let $E_4$ denotes $4\times4$-unit matrix, and let  $C(\nu)$
is result of replacement of first row of the matrix $E_4$
 by the row in (\ref{eq:113bd1}).
 Let further $D(\nu)$  denotes
the adjoint matrix to the matrix $C(\nu).$
 Then
\begin{equation}\label{eq:15.1.b}
C(\nu)=\left(\begin{matrix}
      r_1(\nu) &     r_2(\nu) &     r_3(\nu) &     r_4(\nu)\\
      0        &     1        &     0        &     0        \\
      0        &     0        &     1        &     0        \\
      0        &     0        &     0        &     1
 \end{matrix}\right),
 \end{equation}
\begin{equation}\label{eq:15.1.c}
D(\nu)=\left(\begin{matrix}
       1&-r_2(\nu)      &-r_3(\nu)&-r_4(\nu)\\
       0& r_1(\nu)      &      0        &      0        \\
       0&     0         & r_1(\nu)      &      0        \\
       0&     0         &      0        &     r_1(\nu)
\end{matrix}\right).
\end{equation}
Clearly,
\begin{equation}\label{eq:15.1.d}
C(\nu)D(\nu)=(\mu(\nu))^2E_4,\,C(-\nu-2)=C(\nu),\,D(-\nu-2)=D(\nu).
\end{equation}
Let
\begin{equation}\label{eq:15.1.e}
A^{\ast\ast}(z;\nu)=C(\nu-1)A^{\ast}_{1,0}(z;\nu)D(\nu).
\end{equation}
Then
\begin{equation}\label{eq:15.1.f}
A^{\ast\ast}(z;-\nu-1)=C(-\nu-2)A^{\ast}(z;-\nu-1)D(-\nu-1)=
\end{equation}
$$C(\nu)A^{\ast}(z;-\nu-1)D(\nu-1),$$
and, in view of (\ref{eq:15.1.d}), (\ref{eq:93cc}), (\ref{eq:15.1.e}),
\begin{equation}\label{eq:15.1.g}
A^{\ast\ast}(z;-\nu-1)A^{\ast\ast}_{1,0}(z;\nu)=
\end{equation}
$$C(\nu)A^{\ast}(z;-\nu-1)D(\nu-1)
C(\nu-1)A^{\ast}(z;\nu)D(\nu)=$$
$$-(\mu(\nu)\mu(\nu-1))^2(\nu(\nu+1))^5E_4.$$
Let
\begin{equation}\label{eq:15.1.h}
Y_k(z;\nu)=C(\nu)X_k(z;\nu),
\end{equation}
where
$k=1,\,2,\,3,\,\vert z\vert>1,\,
\nu\in M_1^{\ast\ast\ast}=((-\infty,-2]\cup[0,+\infty))\cap{\mathbb Z}.$
Then, in view of (\ref{eq:93be0}), (\ref{eq:15.1.d}), (\ref{eq:93ca}),
\begin{equation}\label{eq:15.1.i}
Y_k(z;-\nu-2)=Y_k(z;\nu),
\end{equation}
\begin{equation}\label{eq:15.1.aj}
A^{\ast\ast}(z;\nu)Y_k(z;\nu)=C(\nu-1)A^{\ast}(z;\nu)D(\nu)C(\nu)X_k(z;\nu)=
\end{equation}
%$$$$
$$\mu(\nu)^2C(\nu-1)A^{\ast}(z;\nu)X_k(z;\nu)=$$
$$\mu(\nu)^2\nu^5C(\nu-1)X_k(z;\nu-1)=\mu(\nu)^2\nu^5Y_k(z;\nu-1),$$
where
$k=1,\,2,\,3,\,\vert z\vert>1,\,
\nu\in M_1^{\ast}=((-\infty,-2]\cup[1,+\infty))\cap{\mathbb Z}.$
Replacing in the equality (\ref{eq:15.1.aj})
$\nu\in M_1^{\ast}=((-\infty,-2]\cup[1,+\infty))\cap{\mathbb Z}$
 by
$$\nu:=-\nu-2\in M_1^{\ast\ast}=((-\infty,-3]\cup[0,+\infty))\cap{\mathbb Z},$$
and taking in account (\ref{eq:15.1.i}) we obtain
the equality
\begin{equation}\label{eq:15.1.ba}
-A^{\ast\ast}(z;-\nu-2)Y_k(z;\nu)=\mu(\nu)^2(\nu+2)^5Y_k(z;\nu+1),
\end{equation}
where $k=1,\,2,\,3,\,\vert z\vert>1,\,
\nu\in M_1^{\ast\ast}=((-\infty,-3]\cup[0,+\infty))\cap{\mathbb Z}.$
%%%%%%%%%%%%%%%%%%%%%%%%%%%%%%%%%%%%%%\S 3%%%%%%%%%%%%%%%%%%%%%%%%%%%%%%%
\refstepcounter{section}
{\begin{center}\large\bf
\S 3. Calculation of the matrix $A^{\ast\ast}_{1,0}(z;\nu).$
\end{center}}
We denote by
%$a_{1,0,i,j}^{\vee\vee}(1;\nu),\,
%a_{1,0,i,j}^{\wedge\wedge}(1;\nu)$ and
$a_{i,j}^{\ast\ast}(1;\nu),$
where $i,j=1,2,3,4,$
the expressions, which stand on intersection of $i$-th row and $j$-th column
in the matrix $A^{ast\ast}(1;\nu).$ Let
\setcounter{equation}0
\begin{equation}\label{eq:15.2.a}
V^{\ast\ast}(\nu)=C(\nu-1)V^{\ast}(\nu)D(\nu).
\end{equation}
Then, in view of $(\ref{eq:15d}),$
 \begin{equation}\label{eq:15.2.b}
A^{\ast\ast}(z;\nu)=
A^{\ast\ast}(1;\nu)+(z-1)V^{\ast\ast}(\nu),
\end{equation}
where the matrix $V^{\ast\ast}(\nu)$ does not depend from $z.$
Clearly, the first row of the matrix $C(\nu-1)A^{\ast}(1,\nu)$
coincides with the row $R(\nu-1)A^{\ast}(z,\nu)$
and, according to the Theorem 2 coincides
 with the row $\nu^5R(\nu),$ i.e. with
the first row of the matrix $\nu^5C(\nu).$
 Therefore, in view
of (\ref{eq:15.1.d}), the first row of the matrix
 $A^{\ast\ast}(1,\nu)$ is equal to $(\mu_1(\nu)^2)\nu^5\bar e_{4,1},$
where $\bar e_{4,l}$ denotes the $l$-th row of the matrix $E_4$
for $l=1,\,2,\,3,\,4.$
Hence
%\begin{equation}\label{eq:15.2.a11}
%a_{1,1}^{\ast\ast}(1;\nu)=\tau^4(\tau-1)^5,
%\end{equation}
\begin{equation}\label{eq:15.2.vw1,234}
a_{1,1}^{\ast\ast}(1;\nu)=\tau^4(\tau-1)^5,\,a_{1,k}^{\ast\ast}(1;\nu)=0,\,
\text{where}\,\, k=2,\,3,\,4.
\end{equation}
%where $k=2,\,3,\,4.$
 Clearly, the second,  third and fourth row of the
  matrix $C(\nu-1)A^{\ast}(1,\nu)$
  coincides with respectively  the second, third and fourth row
 of  $A^{\ast}(1,\nu).$

In view of (\ref{eq:a21a}), (\ref{eq:a31a}) and (\ref{eq:a41a}),
\begin{equation}\label{eq:a21asas}%!
a_{2,1}^{\ast\ast}(1;\nu)=a_{2,1}^\ast(1;\nu)=
-(12\tau^6-26\tau^5+20\tau^4-7\tau^3+\tau^2)=-\tau^2\times
\end{equation}
$$\,\,\,\,\,\,\,\,(\tau-1)(12\tau^3-14\tau^2+6\tau-1)=
-\tau^2(\tau-1)(2\tau-1)(6\tau^2-4\tau+1).$$
\begin{equation}\label{eq:a31asas}%!
a_{3,1}^{\ast\ast}(1;\nu)=a_{3,1}^\ast(1;\nu)=
8\tau^7-26\tau^6+33\tau^5-21\tau^4+7\tau^3-\tau^2=\tau^2\times
\end{equation}
$$(\tau-1)(8\tau^4-18\tau^3+15\tau^2-6\tau+1)=\tau^2(\tau-1)^2\times$$
$$(8\tau^3-10\tau^2+5\tau-1)=\tau^2(\tau-1)^2(2\tau-1)(4\tau^2-3\tau+1).$$
\begin{equation}\label{eq:a41asas}%!
a_{4,1}^{\ast\ast}(1;\nu)=-\tau^2(\tau-1)^3(2\tau-1)(2\tau^2-2\tau+1)=
\end{equation}
$$a_{4,1}^\ast(1;\nu)=
-\tau^2(\tau-1)(4\tau^5-14\tau^4+20\tau^3-15\tau^2+6\tau-1)=$$
$$-\tau^2(\tau-1)^2(4\tau^4-10\tau^3+10\tau^2-5\tau+1)=$$
$$-\tau^2(\tau-1)^3(4\tau^3-6\tau^2+4\tau-1).$$
In view of (\ref{eq:15.1.b}), (\ref{eq:15.1.c}),  (\ref{eq:15.1.e})
\begin{equation}\label{eq:av234,2}%!!
a_{k,j}^{\ast\ast}(1;\nu)=-r_j(\nu)a_{k,1}^\ast(1;\nu)+
r_{1}(\nu)a_{k,j}^\ast(1;\nu)
\end{equation}
%$$=\mu^2a_{1,0,k,2+2j}^\ast(1;\nu),$$
where $j,\,k=2,\,3,\,4.$

%In view of (\ref{eq:a12a}),
%\begin{equation}\label{eq:a12az1}%!
%a_{1,2}^\ast(1;\nu)=-20\tau^2-24\tau^4+4\tau+36\tau^3.
%\end{equation}

In view of (\ref{eq:a22a}), (\ref{eq:a32a}), (\ref{eq:a42a}),
  (\ref{eq:a24a}),  (\ref{eq:a34a}),
  (\ref{eq:av234,2}) and (\ref{eq:93bd}),
\begin{equation}\label{eq:a22asas}%!
a_{2,2}^{\ast\ast}(1;\nu)=\tau^4a_{2,2}^\ast(1;\nu)=
\tau^5(\tau-1)(17\tau^3-24\tau^2+12\tau-2)=
\end{equation}
$$(17\tau^5-41\tau^4+36\tau^3-14\tau^2+2\tau)\tau^4.$$

\begin{equation}\label{eq:a32asas}%!
a_{3,2}^{\ast\ast}(1;\nu)=\tau^4a_{3,2}^\ast(1;\nu)=
-2\tau^5(\tau-1)^2(6\tau^3-9\tau^2+5\tau^2-1)=
\end{equation}
$$-(12\tau^6+42\tau^5-58\tau^4+40\tau^3-14\tau^2+2\tau)\tau^4=$$
$$-2\tau^5(6\tau^5-21\tau^4+29\tau^3-20\tau^2+7\tau-1)=$$
$$-2\tau^5(\tau-1)(6\tau^4-15\tau^3+14\tau^2-6\tau+1),$$

\begin{equation}\label{eq:a42asas}%!
a_{4,2}^{\ast\ast}(1;\nu)=\tau^4a_{4,2}^{\ast}(1;\nu)=
\tau^5(\tau-1)^3(8\tau^3-14\tau^2+9\tau-2)=
\end{equation}
$$(8\tau^7-38\tau^6+75\tau^5-79\tau^4+47\tau^3-15\tau^2+2\tau)\tau^4=$$
$$\tau^5(\tau-1)(8\tau^5-30\tau^4+45\tau^3-34\tau^2+13\tau-2)=$$
$$\tau^5(\tau-1)^2(8\tau^4-22\tau^3+23\tau^2-11\tau+2),$$
%$$,$$
\begin{equation}\label{eq:1abk4}%!
a_{k,4}^{\ast\ast}(1;\nu)=\tau^4a_{1,0,k,4}^\ast(1;\nu)=0
\end{equation}
for $k=2,\,3.$ In view of (\ref{eq:a44a}), (\ref{eq:av234,2})
 and (\ref{eq:93bd}),
\begin{equation}\label{eq:a44az1}%!
a_{4,4}^{\ast\ast}(1;\nu)=\tau^4a_{4,4}^\ast(1;\nu)=\tau^4\times
\end{equation}
$$(\tau^5-5\tau^4+10\tau^3-10\tau^2+5\tau-1)=\tau^4(\tau-1)^5.$$
In view of (\ref{eq:a23a}),
\begin{equation}\label{eq:a23az1}%!
a_{2,3}^\ast(1;\nu)=8\tau^3-14\tau^2+7\tau-1=(\tau-1)(8\tau^2-6\tau+1),
\end{equation}

In view of (\ref{eq:a33a}),
\begin{equation}\label{eq:a33az1}%1
a_{\alpha,0,3,3}^\ast(1;\nu)=
 \tau^5-9\tau^4+20\tau^3-18\tau^2+7\tau-1=
\end{equation}
$$(\tau-1)(\tau^4-8\tau^3+12\tau^2-6\tau+1=$$
$$(\tau-1)^2(\tau^3-7\tau^2+5\tau-1),$$%%
In view of (\ref{eq:a43a}),
\begin{equation}\label{eq:a43az1}%!
a_{4,3}^\ast(1;\nu)=-4\tau^6+22\tau^5-48\tau^4+53\tau^3-31\tau^2+9\tau-1=
\end{equation}
$$(\tau-1)(-4\tau^5+18\tau^4-30\tau^3+23\tau^2-8\tau+1)=$$
$$(\tau-1)^2(-4\tau^4+14\tau^3-16\tau^2+7\tau-1)=$$
$$(\tau-1)^3(-4\tau^3+10\tau^2-6\tau+1)$$
In view of (\ref{eq:a21asas}), (\ref{eq:a23az1}),
 (\ref{eq:a31asas}), (\ref{eq:a33az1}),
(\ref{eq:a41asas}), (\ref{eq:a43az1}),
 (\ref{eq:av234,2})  and (\ref{eq:93bd}),
\begin{equation}\label{eq:a23asas}
a_{2,3}^{\ast\ast}(1;\nu)=2\tau^2a_{2,1}^{\ast}(1;\nu)+
\tau^4a_{2,3}^{\ast}(1;\nu)=(\tau-1)\times
\end{equation}
$$
(-2\tau^4(2\tau-1)(6\tau^2-4\tau+1)+\tau^4(2\tau-1)(4\tau-1))=
$$
$$\tau^4(\tau-1)(2\tau-1)(-12\tau^2+12\tau-3)=-3\tau^4(2\tau-1)^3,$$
\begin{equation}\label{eq:a33asas}
a_{3,3}^{\ast\ast}(1;\nu)=2\tau^2a_{3,1}^{\ast}(1;\nu)+
\tau^4a_{3,3}^{\ast}(1;\nu)=(\tau-1)^2\times
\end{equation}
$$(2\tau^4(8\tau^3-10\tau^2+5\tau-1)+\tau^4(\tau^3-7\tau^2+5\tau-1))=$$
$$(\tau-1)^2(17\tau^3-27\tau^2+15\tau-3)=
\tau^4(\tau-1)^2((\tau-1)^3+2(2\tau-1)^3),$$
\begin{equation}\label{eq:a43asas}
a_{4,3}^{\ast\ast}(1;\nu)=2\tau^2a_{4,1}^{\ast}(1;\nu)+
\tau^4a_{4,3}^{\ast}(1;\nu)=(\tau-1)^3\times
\end{equation}
$$(-2\tau^4(4\tau^3-6\tau^2+4\tau-1)+\tau^4(-4\tau^3+10\tau^2-6\tau+1))=$$
$$-\tau^4(\tau-1)^3(12\tau^3-22\tau^2+14\tau-3)=$$
%(\tau-1/2)(12\tau^2-16\tau+6)
$$-\tau^4(\tau-1)^3(2\tau-1)(6\tau^2-8\tau+3).$$
%\begin{equation}\label{eq:1as43}
%%%%%%%%%%%%%%%%%%%%%%%%%%%%%%%%%%%%%%%\S 4%%%%%%%%%%%%%%%%%%%%%%%%%%%%%%%%%%%
\refstepcounter{section}
{\begin{center}\large\bf
\S 4. Properties of the functions considered in \S 1.\end{center}}
The function $t^r(R(t,\nu))^2$ (see (\ref{eq:d}))
 is regular at $t=\infty$ for $r=0,\,1,\,2,$ and has a pole of first order
 at $t=\infty$ for $r=3.$ So, in the case $r=0,\,1,\,2$ we
have the equalities
\setcounter{equation}0
\begin{equation}\label{eq:15.4.a.0a}
\Res(t^r(R(t,\nu))^2,t=\infty)=-[r/3]\,\,\text{for}\,\,r=0,\,1,\,2,\,3,
\end{equation}
\begin{equation}\label{eq:15.4.a.0b}
\lim\limits_{t\to \infty} t^r(R(t,\nu))^2=0\,\,\text{for}\,\,r=0,\,1,\,2,\,3.
\end{equation}
%and in the case $r=3,\,\alpha=1$ we
%have the equalities
%\begin{equation}\label{eq:15.4.a.1a}
%\Res(t^3(R(1,t,\nu))^2,t=\infty)=-1,\,
%\end{equation}
%\begin{equation}\label{eq:15.4.a.1b}
%\lim\limits_{t\to \infty} t^3(R(1,t,\nu))^2=0.
%\end{equation}
 In view of (\ref{eq:h}),
\begin{equation}\label{eq:15.4.hr}
\delta^rf_2^{\ast}(z, \nu)= \sum\limits_{t=1}^{+\infty}z^{-t}
(\nu+1)^2(-t)^r(R(t,\nu))^2,
\end{equation}
where we consider $r=0,\,1,\,2,\,3.$
Expanding  $(\nu+1)^2(-t)^r(R(t,\nu))^2$
 into partial fractions relatively $t$, we obtain
\begin{equation}\label{eq:6bg}
(\nu+1)^2(-t)^r(R(t;\nu))^2=\sum\limits_{i=1}^2
\left(\sum\limits_{k=0}^{\nu+1}
\beta_{i,k,\nu}^{(r)}(t + k)^{-i}\right),
\end{equation}
%$$\sum\limits_{i=1}^2
%\left(\sum\limits_{k=0}^{\nu+1}
%\beta_{i,k,\nu}^{(r)}(t + k)^{-i}\right),$$
where $\nu\in{\mathbb N}_0,\,r=0,\,1,\,2,\,3,$
\begin{equation}\label{eq:6bh}
\beta_{2-j,k,\nu}^{(r)}=(\nu+1)^2%\times
\frac1{j!}\lim\limits_{t\to-k}\left(\frac{\partial}{\partial t}\right)^j
((-t)^r(R(t,\nu)(t+k))^2)
\end{equation}
%$$$$
 for $j=0,\,1.$
In view of (\ref{eq:15.4.a.0a}) and (\ref{eq:6bg}),
\begin{equation}\label{eq:6bh.0}
\sum\limits_{k=0}^{\nu+1}\beta_{1,k,\nu}^{(r)}=-[r/3](\nu+1)^2
\,\,\text{for}\,\, r=0,\,1,\,2,\,3.
\end{equation}
% for $r=0,\,1,\,2,$
%\begin{equation}\label{eq:6bh.1}
%\sum\limits_{k=0}^{\nu+1}\beta_{1,k,\nu}^{(3)}=-(\nu+1)^2.
%\end{equation}
In view of (\ref{eq:6bg}),
\begin{equation}\label{eq:6bgf}
-(\nu+1)^2\frac\partial
{\partial t}((-t)^r(R(t;\nu))^2)=%$$
\sum\limits_{i=1}^2
\left(\sum\limits_{k=0}^{\nu+1}
\beta_{i,k,\nu}^{(r)}i(t + k)^{-i-1}\right),%$$
\end{equation}
where $\nu\in{\mathbb N}_0,\,r=0,\,1,\,2,\,3.$
 Let
\begin{equation}\label{eq:6bi0}
S_{i,k}(\nu)=-\left(\sum\limits_{\kappa=k+1}^{\nu+k}1/\kappa^i\right)-%$$
\left(\sum\limits_{\kappa=1}^{\nu+1-k}1/\kappa^i\right)+
\sum\limits_{\kappa=1}^k1/\kappa^i,%$$
\end{equation}
where
$\nu\in{\mathbb N}_0,\,i\in{\mathbb N},\,k\in[0,\,\nu+1]\cap{\mathbb Z}.$
In particular,
%\begin{equation}\label{eq:6bi1a0}
%S_{1,0}(0,0)=0,\,S_{1,0}(0,1)=-2,\,S_{1,1}(0,1)=1/2,
%\end{equation}
\begin{equation}\label{eq:6bi1k01a1n0}
S_{1,0}(0)=-1,\,S_{1,1}(0)=1,
\end{equation}
\begin{equation}\label{eq:6bi1k012a1n1}
S_{1,0}(1)=-5/2,\,S_{1,1}(1)=-1/2,\,S_{1,2}(1)=7/6,
\end{equation}
\begin{equation}\label{eq:6bi1k0a1n2}
S_{1,0}(2)=-(1+1/2)-(1+1/2+1/3)=-10/3,
\end{equation}
\begin{equation}\label{eq:6bi1k1a1n2}
S_{1,1}(2)=-(1/2+1/3)-(1+1/2)+1=-4/3
\end{equation}
\begin{equation}\label{eq:6bi1k2a1n2}
S_{1,2}(2)=-(1/3+1/4)-1+(1+1/2)=-1/12,
\end{equation}
\begin{equation}\label{eq:6bi1k3a1n2}
S_{1,3}(2)=-(1/4+1/5)+(1+1/2+1/3)=83/60.
\end{equation}
In view of (\ref{eq:6bh}), (\ref{eq:d}) and (\ref{eq:6bi0})
\begin{equation}\label{eq:6bi2}
\beta_{2,k,\nu}^{(0)}=
 \left(\frac{(\nu+k)!}{k!}\times\frac{\nu+1}{(\nu+1-k)!k!}\right)^2=
\binom{\nu+1}k^2\binom{\nu+k}k^2,
\end{equation}
\begin{equation}\label{eq:6bi1}
\beta_{1,k,\nu}^{(0)}=2\beta_{2,k,\nu}^{(0)}S_{1,k}(\nu),
\end{equation}
where
$\nu\in{\mathbb N}_0,\,i\in{\mathbb N},\,k\in[0,\,\nu+1]\cap{\mathbb Z}.$
In particular,
\begin{equation}\label{eq:r0a1i2k01n0}
\beta_{2,0,0}^{(0)}=\beta_{2,1,0}^{(0)}=\beta_{2,0,1}^{(0)}=1,
\,\beta_{2,1,1}^{(0)}=16,\,\beta_{2,2,1}^{(0)}=9,
\end{equation}
\begin{equation}\label{eq:r0a1i2k012n1}
\beta_{2,0,1}^{(0)}=1,\,\beta_{2,1,1}^{(0)}=16,\,\beta_{2,2,1}^{(0)}=9,
\end{equation}
\begin{equation}\label{eq:r0a1i2k01n2}
\beta_{2,0,2}^{(0)}=1,\,\beta_{2,1,2}^{(0)}=81,\,
\beta_{2,2,2}^{(0)}=324,\,\beta_{2,3,2}^{(0)}=100.
\end{equation}
\begin{equation}\label{eq:r0a1i2k23n2}
\beta_{2,2,2}^{(0)}=324,\,\beta_{2,3,2}^{(0)}=100.
\end{equation}
In view of (\ref{eq:6bi1}), (\ref{eq:r0a1i2k01n0}) -- (\ref{eq:r0a1i2k01n2})
 and (\ref{eq:6bi1k01a1n0}) -- (\ref{eq:6bi1k3a1n2}),
\begin{equation}\label{eq:r0a1i1k0n0}
\beta_{1,k,0}^{(0)}=2\beta_{2,k,0}^{(0)}S_{1,k}(0)=-2\times1\times(-1)^k\,\,
\text{for}\,\,k=0,\,1,
\end{equation}
\begin{equation}\label{eq:r0a1i1k0n1}
\beta_{1,0,1}^{(0)}=
2\beta_{2,0,1}^{(0)}S_{1,0}(1)=2\times1\times(-5/2)=-5,
\end{equation}
\begin{equation}\label{eq:r0a1i1k1n1}
\beta_{1,1,1}^{(0)}=2\beta_{2,1,1}^{(0)}S_{1,1}(1)=
2\times16\times(-1/2)=-16,
\end{equation}
\begin{equation}\label{eq:r0a1i1k2n1}
\beta_{1,2,1}^{(0)}=2\beta_{2,2,1}^{(0)}S_{1,2}(1)=2\times9\times(7/6)=21,
\end{equation}
\begin{equation}\label{eq:r0a1i1k0n2}
\beta_{1,0,2}^{(0)}=2\beta_{2,0,2}^{(0)}S_{1,0}(2)=2(-10/3)=-20/3,
\end{equation}
\begin{equation}\label{eq:r0a1i1k1n2}
\beta_{1,1,2}^{(0)}=2\beta_{2,1,2}^{(0)}S_{1,1}(2)= 2\times81\times(-4/3)=-216,
\end{equation}
\begin{equation}\label{eq:r0a1i1k2n2}
\beta_{1,2,2}^{(0)}=2\beta_{2,2,2}^{(0)}S_{1,2}(2)=
2\times324\times(-1/12)=-54,
\end{equation}
\begin{equation}\label{eq:r0a1i1k3n2}
\beta_{1,3,2}^{(0)}=2\beta_{2,3,2}^{(0)}S_{1,3}(2)=
2\times100\times(83/60)=830/3.
\end{equation}
We put in (\ref{eq:6bg}) $r=0,$ and multiply both sides of obtained
equality by $(-t)^r$ for $r=0,\,1,\,2,\,3.$ Then we see that
\begin{equation}\label{eq:6bg.1}
-t(\nu+1)^2(R(t;\nu))^2=\sum\limits_{i=1}^2
\left(\sum\limits_{k=0}^{\nu+1} \frac{\beta_{i,k,\nu}^{(0)}(-t-k+k)}{(t +
k)^i}\right)=
\end{equation}
$$ \left(\sum\limits_{k=0}^{\nu+1} \frac{k\beta_{2,k,\nu}^{(0)}}{(t
+ k)^2}\right)+
\left(\sum\limits_{k=0}^{\nu+1}
\frac{k\beta_{1,k,\nu}^{(0)}-\beta_{2,k,\nu}^{(0)}}{t+k}\right)-
\sum\limits_{k=0}^{\nu+1} \beta_{1,k,\nu}^{(0)},$$
\begin{equation}\label{eq:6bg.2}
(-t)^2(\nu+1)^2(R(\alpha,t;\nu))^2=\sum\limits_{i=1}^2
\left(\sum\limits_{k=0}^{\nu+1}
\frac{\beta_{i,k,\nu}^{(0)}(t+k-k)^2}{(t+k)^i}\right)=
\end{equation}
$$
\left(\sum\limits_{k=0}^{\nu+1}
\frac{k^2\beta_{2,k,\nu}^{(0)}}{(t + k)^2}\right)+
\left(\sum\limits_{k=0}^{\nu+1}\frac{
k^2\beta_{1,k,\nu}^{(0)}-2k\beta_{2,k,\nu}^{(0)}}
{t+k}\right)+
\sum\limits_{k=0}^{\nu+1}(
\beta_{2,k,\nu}^{(0)}+(t-k)\beta_{1,k,\nu}^{(0)}),$$
\begin{equation}\label{eq:6bg.3}
(-t)^3(\nu+1)^2(R(\alpha,t;\nu))^2=\sum\limits_{i=1}^2
\left(\sum\limits_{k=0}^{\nu+1}
\frac{\beta_{i,k,\nu}^{(0)}(-t-k+k)^3}{(t+k)^i}\right)=
\end{equation}
$$
\left(\sum\limits_{k=0}^{\nu+1}\frac
{k^3\beta_{2,k,\nu}^{(0)}}{(t + k)^2}\right)+
\left(\sum\limits_{k=0}^{\nu+\alpha}\frac
{k^3\beta_{1,k,\nu}^{(0)}-3k^2\beta_{2,k,\nu}^{(0)}}
{t + k}\right)-$$
$$\left(\sum\limits_{k=0}^{\nu+1}
(t-2k)(\beta_{2,k,\nu}^{(0)}\right)-
\sum\limits_{k=0}^{\nu+1}
(t^2-kt+k^2)\beta_{1,k,\nu}^{(0)}.$$
The equality (\ref{eq:6bh.0}) with $r=0$ again follows
 from (\ref{eq:15.4.a.0b}) with $r=1$ and (\ref{eq:6bg.1}); moreover,
in view of (\ref{eq:6bg}) with $r=1,$ and (\ref{eq:6bg.1}),
%\begin{equation}\label{eq:6bg.1.2}
%\beta_{2,k,\nu}^{(1)}=k\beta_{2,k,\nu}^{(0)},\,
%\end{equation}
\begin{equation}\label{eq:6bg.1.1}
\beta_{2,k,\nu}^{(1)}=k\beta_{2,k,\nu}^{(0)},\,
\beta_{1,k,\nu}^{(1)}=k\beta_{1,k,\nu}^{(0)}-
\beta_{2,k,\nu}^{(0)}
\end{equation}
for $k=0,\,...,\,\nu+1.$
The equality (\ref{eq:6bh.0}) with $r=1$ again follows
 from (\ref{eq:15.4.a.0b}) with $r=2,$ (\ref{eq:6bh.0}) with $r=0,$
(\ref{eq:6bg.2}) and (\ref{eq:6bg.1.1}); moreover,
in view of (\ref{eq:6bg}) with $r=2,$ and (\ref{eq:6bg.2}),
%\begin{equation}\label{eq:6bg.2.2}
%\beta_{2,k,\nu}^{(2)}=k^2\beta_{2,k,\nu}^{(0)},\,
%\end{equation}
\begin{equation}\label{eq:6bg.2.1}
\beta_{2,k,\nu}^{(2)}=k^2\beta_{2,k,\nu}^{(0)},\,
\beta_{1,k,\nu}^{(2)}=k^2\beta_{1,k,\nu}^{(0)}-2k\beta_{2,k,\nu}^{(0)}
\end{equation}
for $k=0,\,...,\,\nu+1.$
The equality (\ref{eq:6bh.0}) with $r=2$ again follows
 from  (\ref{eq:15.4.a.0b}), (\ref{eq:6bh.0}) with both $r\in\{0,\,1\},$
(\ref{eq:6bg.3}), (\ref{eq:6bg.1.1}) and from (\ref{eq:6bg.2.1}); moreover,
in view of (\ref{eq:6bg}) with $r=3,$ and (\ref{eq:6bg.3}),
%\begin{equation}\label{eq:6bg.3.2}
%\beta_{2,k,\nu}^{(3)}=k^3\beta_{2,k,\nu}^{(0)},\,
%\end{equation}
\begin{equation}\label{eq:6bg.3.1}
\beta_{2,k,\nu}^{(3)}=k^3\beta_{2,k,\nu}^{(0)},\,
\beta_{1,k,\nu}^{(3)}=k^3\beta_{1,k,\nu}^{(0)}-3k^2\beta_{2,k,\nu}^{(0)}
\end{equation}
for $k=0,\,...,\,\nu+1.$
In view of (\ref{eq:h}) -- (\ref{eq:ac}),
\begin{equation}\label{eq:acrf}
(\delta^r)f_3^\ast(z,\nu)=
(\log(z))(\delta)^r)f_2^{\ast}(z,\nu)+
\end{equation}
$$r(\delta)^{r-1}f_2^{\ast}(z,\nu)
+(\delta)^rf_4^{\ast}(z,\nu)=
(\log(z))(\delta)^r)f_2^{\ast}(z,\nu)+$$
$$\sum\limits_{t=1}^{+\infty} z^{-t}(\nu+1)^2
\left (r(-t)^{r-1}-(-t)^r\frac\partial{\partial t}\right)R^2(t,\nu)=$$
$$(\log(z))(\delta)^r)f_2^{\ast}(z,\nu)-
\sum\limits_{t=1}^{+\infty} z^{-t}
(\nu+1)^2\frac\partial{\partial t}((-t)^rR^2)(t,\nu).$$
In view of (\ref{eq:6bg}), (\ref{eq:6bgf}), (\ref{eq:15.4.hr}) and
(\ref{eq:acrf}),
\begin{equation}\label{eq:15.4.hr1f}
\delta^rf_{2+j}^\ast(z;\nu)-
j(\log(z))\delta^rf_2^{\ast}(z;\nu)=
\end{equation}
$$\sum\limits_{i=1}^2
\left(\sum\limits_{t=1}^{+\infty}
\left(\sum\limits_{k=0}^{\nu+1}(1-j+ij)\beta_{i,k,\nu}^{(r)}
z^kz^{-t-k}(t + k)^{-i-j}\right)\right)=
$$
$$\sum\limits_{i=1}^2\left(\sum\limits_{k=0}^{\nu+1}
(1-j+ij)\beta_{i,k,\nu}^{(r)}z^k
\left(\sum\limits_{t=1}^{+\infty}
z^{-t-k} (t + k)^{-i-j}\right)\right)=
$$
$$\sum\limits_{i=1}^2\left(\sum\limits_{k=0}^\nu
(1-j+ij)\beta_{i,k,\nu}^{(r)}z^k\left(L_{i+1}(1/z)-
\sum\limits_{\tau=1}^kz^{-\tau}i(\tau)^{-i-j}\right)\right)=
$$
$$\left(\sum\limits_{i=1}^2
(1-j+ij)\beta^{\ast(r)}_i(z;\nu)L_{i+j}(1/z)\right)-
\beta^{\ast(r)}_{3+j}(z;\nu),$$
where $j=0,\,1,\,r=0,\,1,\,2,\,3,\,\vert z\vert>1,$
\begin{equation}\label{eq:6cd}
L_s(1/z)=\sum\limits_{n=1}^\infty 1/(z^nn^s),\,
\beta^{\ast(r)}_i(z;\nu)=
\sum\limits_{k=0}^{\nu+1}\beta_{i,k,\nu}^{(r)}z^k,
\end{equation}
for $s\in{\mathbb Z},\,i\in\{1,2\},\,\nu\in{\mathbb N}_0,$
\begin{equation}\label{eq:6ce}
\beta^{\ast(r)}_{3+j}(z;\nu)=\sum\limits_{i=1}^2\left(\sum\limits_{k=0}^{\nu+1}
(1-j+ij)\beta_{i,k,\nu}^{(r)}\left(\sum\limits_{\tau=1}^k
z^{k-\tau}(\tau)^{-i-j}\right)\right)=
\end{equation}
%$$\sum\limits_{i=1}^2\left(\sum\limits_{k=0}^{\nu+1}
%(1-j+ij)\beta_{i,k,\nu}^{(r)}
%\left(\sum\limits_{\tau=1}^k
%z^{k-\tau}(\tau)^{-i-j}\right)\right)=
%$$
$$\sum\limits_{\sigma=0}^{\nu+\alpha-1}z^\sigma
\sum\limits_{\tau=1}^{\nu+1-\sigma}
\sum\limits_{i=1}^2(1-j+ij)\beta_{i,\sigma+\tau,\nu}^{(r)}
(\tau)^{-i-j}.$$
In fiew of (\ref{eq:b}),
\begin{equation}\label{eq:b1}
f_1^{\ast}(z, \nu)=\beta^{\ast(0)}_2(z;\nu).
\end{equation}
In view of (\ref{eq:6bh.0}) and (\ref{eq:6cd}), if
 $r=0,\,1,\,2,$ then
\begin{equation}\label{eq:6cd1}
\beta^{\ast(r)}_1(z;\nu)=
(z-1)\beta^{\ast\vee(r)}_1(z;\nu),
\end{equation}
where $\beta^{\ast\vee(r)}_1(z;\nu)\in{\mathbb Q}[z],$
when $\nu\in{\mathbb N}_0.$
In view of (\ref{eq:6bh.0}) and (\ref{eq:6cd}),
\begin{equation}\label{eq:6cd2}
\beta^{\ast(3)}_1(z;\nu)=-(\nu+1)^2+
(z-1)\beta^{\ast\vee(3)}_1(z;\nu),
\end{equation}
where $\beta^{\ast\vee(3)}_1(z;\nu)\in{\mathbb Q}[z],$
when $\nu\in{\mathbb N}_0.$
In view of (\ref{eq:6bg.1.1}) -- (\ref{eq:6bg.3.1}), (\ref{eq:6cd}),
\begin{equation}\label{eq:6bg.1.2a}
\beta_2^{\ast(1)}(z;\nu)=\delta\beta_2^{\ast(0)}(z;\nu)=
\delta f_1^{\ast}(z, \nu),\,
\beta_1^{\ast(1)}(z;\nu)=
\end{equation}
$$\delta\beta_1^{\ast(0)}(z;\nu)-
\beta_2^{\ast(0)}(z;\nu),$$
%\begin{equation}\label{eq:6bg.1.1a}
%\beta_1^{\ast(1)}(z;\nu)=\delta\beta_1^{\ast(0)}(z;\nu)-
%\beta_2^{\ast(0)}(z;\nu),
%\end{equation}
%\begin{equation}\label{eq:6bg.2.2a}
%\beta_2^{\ast(2)}(z;\nu)=
%\delta^2\beta_2^{\ast(0)}(z;\nu),\,
%\end{equation}
\begin{equation}\label{eq:6bg.2.1a}
\beta_2^{\ast(2)}(z;\nu)=
\delta^2\beta_2^{\ast(0)}(z;\nu)=\delta^2f_1^{\ast}(z, \nu),\,
\beta_1^{\ast(2)}(z;\nu)=
\end{equation}
$$\delta^2\beta_1^{\ast(0)}(z;\nu)-
2\delta\beta_2^{\ast(0)}(z;\nu),$$
%\begin{equation}\label{eq:6bg.3.2a}
%\beta_2^{\ast(3)}(z;\nu)=
%\delta^3\beta_2^{\ast(0)}(z;\nu),\,
%\end{equation}
\begin{equation}\label{eq:6bg.3.1a}
\beta_2^{\ast(3)}(z;\nu)=
\delta^3\beta_2^{\ast(0)}(z;\nu),\,\beta_1^{\ast(3)}(z;\nu)=
\delta^3\beta_1^{\ast(0)}(z;\nu)-
3\delta^2\beta_2^{\ast(0)}(z;\nu).
\end{equation}
 Clearly,
\begin{equation}\label{eq:6cg}
(-\delta)^kL_n(1/z)=L_{n-k}(1/z),
\end{equation}
where $k\in[0,+\infty)\cap{\mathbb Z},\,n\in{\mathbb Z},\,\vert z\vert>1,$
\begin{equation}\label{eq:6ch0}
L_1(1/z)=-\log(1-1/z),\,-\delta L_1(1/z)=1/(z-1)=
\end{equation}
$$L_0(1/z),\,\delta^2L_1(1/z)=1/(z-1)+1/(z-1)^2=$$
$$L_{-1}(1/z),\,-\delta^3L_1(1/z)=L_{-2}(1/z)=1/(z-1)+3/(z-1)^2+2/(z-1)^3.$$
We apply the operator $\delta$ to the equality (\ref{eq:15.4.hr1f})
 for $r=0,\,1,\,2.$ Then,
 in view of  (\ref{eq:6cg}), we obtain the equality
\begin{equation}\label{eq:15.4.hr2f}
\delta^{r+1}f_{2+j}^{\ast}(z;\nu)-
j(\log(z))\delta^{r+1}f_2^{\ast}(z, \nu)=j\delta^rf_2^{\ast}(z;\nu)+
\end{equation}
%$$j\delta^rf_2^{\ast}(z;\nu)+$$
$$\left(\sum\limits_{i=1}^2
((1-j+ij)\delta\beta^{\ast(r)}_i(z;\nu))L_{i+j}(1/z)\right)-
\delta\beta^{\ast(r)}_{3+j}(z;\nu)-$$
$$\left(\sum\limits_{i=1}^2
(1-j+ij)\beta^{\ast(r)}_i(z;\nu)L_{i+j-1}(1/z)\right)=$$
$$j\left(\left(\sum\limits_{i=1}^2
\beta^{\ast(r)}_i(z;\nu)L_i(1/z)\right)-
\beta^{\ast(r)}_3(z;\nu)\right)+$$
$$\left(\sum\limits_{i=1}^2
((1-j+ij)\delta\beta^{\ast(r)}_i(z;\nu))L_{i+j}(1/z)\right)-
\delta\beta^{\ast(r)}_{3+j}(z;\nu)-$$
$$\left(\sum\limits_{i=1}^2
(1-j+ij)\beta^{\ast(r)}_i(z;\nu)L_{i+j-1}(1/z)\right).$$
It follws from (\ref{eq:15.4.hr2f}) with $j=0$ that
\begin{equation}\label{eq:15.4.hr2f0}
\delta^{r+1}f_2^{\ast}(z;\nu)=-\delta\beta^{(r)}_3(z;\nu)+
\end{equation}
$$\left(\sum\limits_{i=1}^2(\delta\beta^{\ast(r)}_i(z;\nu))L_{i}(1/z)-
\beta^{\ast(r)}_i(z;\nu)L_{i-1}(1/z)\right)=$$
$$(\delta\beta^{\ast(r)}_2(z;\nu))L_2(1/z)+
(\delta\beta^{\ast(r)}_1(z;\nu)-
\beta^{\ast(r)}_2(z;\nu))L_{1}(1/z)-$$
$$\delta\beta^{\ast(r)}_3(z;\nu)-
\beta^{\ast(r)}_1(z;\nu)L_0(1/z).$$
In view of (\ref{eq:15.4.hr1f}) with $j=0,$ (\ref{eq:15.4.hr2f0}),
 (\ref{eq:6cd1}),
\begin{equation}\label{eq:15.4.hr3}
\beta^{\ast(r)}_2(z;\nu)=
\delta\beta^{\ast(r-1)}_2(z;\nu))=
\delta^r\beta^{\ast(0)}_2(z;\nu),
\end{equation}
\begin{equation}\label{eq:15.4.hr4}
\beta^{\ast(r)}_1(z;\nu))=
\delta\beta^{\ast(r-1)}_1(z;\nu)-
\beta^{\ast(r-1)}_2(z;\nu)=
\end{equation}
$$\delta^r\beta^{\ast(0)}_1(z;\nu)-r\delta^{r-1}\beta^{\ast(0)}_2(z;\nu),$$
\begin{equation}\label{eq:15.4.hr5}
\beta^{\ast(r)}_3(z;\nu)=\delta\beta^{\ast(r-1)}_3(z;\nu)+
\beta^{\ast(r-1)}_1(z;\nu))L_0(1/z)=
\end{equation}
$$\delta\beta^{\ast(r-1)}_3(z;\nu)+\beta^{\ast\vee(r-1)}_1(z;\nu),$$
where $r=1,\,2,\,3.$
The equalities (\ref{eq:6bg.1.2a}) -- (\ref{eq:6bg.3.1a})
follow from the equalities (\ref{eq:15.4.hr3}) and (\ref{eq:15.4.hr4}) again.
In view of (\ref{eq:6bi2}), (\ref{eq:6cd}), (\ref{eq:b})
 and (\ref{eq:15.4.hr3})
\begin{equation}\label{eq:15.4.hr8}
\beta^{\ast(r)}_2(z;\nu)=
\delta^rf^{\ast}_1(z;\nu)\in{\mathbb N}[z],
\end{equation}
where $\nu\in{\mathbb N}_0,\,r=0,\,1,\,2,\,3.$
It follws from (\ref{eq:15.4.hr2f}) with $j=1$ that
\begin{equation}\label{eq:15.4.hr2f1}
\delta^{r+1}f_3^{\ast}(z, \nu)=
(\log(z))\delta^{r+1}f_2^{\ast}(z, \nu)+
\end{equation}
$$\left(\left(\sum\limits_{i=1}^2
\beta^{\ast(r)}_i(z;\nu)L_i(1/z)\right)-
\beta^{\ast(r)}_3(z;\nu)\right)+$$
$$\left(\sum\limits_{i=1}^2
i(\delta\beta^{\ast(r)}_i(z;\nu))L_{i+1}(1/z)\right)-
\delta\beta^{\ast(r)}_4(z;\nu)-$$
$$\left(\sum\limits_{i=1}^2
i\beta^{\ast(r)}_i(z;\nu)L_{i}(1/z)\right)=
(\log(z))\delta^{r+1}f_2^{\ast}(z;\nu)+$$
$$\left(\sum\limits_{i=1}^2
i(\delta\beta^{\ast(r)}_i(z;\nu))L_{i+1}(1/z)\right)-
\delta\beta^{\ast(r)}_4(z;\nu)-$$
$$\beta^{\ast(r)}_2(z;\nu)L_{2}(1/z)-
(\beta^{\ast(r)}_3(z;\nu)+
\delta\beta^{(r)}_4(z;\nu))=$$
$$(\log(z))\delta^{r+1}f_2^{\ast}(z;\nu)+
2(\delta\beta^{\ast(r)}_2(z;\nu))L_3(1/z)+$$
$$(\delta\beta^{\ast(r)}_1(z;\nu)-
\beta^{\ast(r)}_2(z;\nu))L_2(1/z)-
(\delta\beta^{\ast(r)}_4(z;\nu)+\beta^{\ast(r)}_3(z;\nu)).$$
In view of (\ref{eq:15.4.hr1f}) with $j=1,$ (\ref{eq:15.4.hr2f1}),
\begin{equation}\label{eq:15.4.hr2f2}
\beta^{\ast(r+1)}_2(z;\nu)=\delta\beta^{\ast(r)}_2(z;\nu)=
\delta^{r+1}\beta^{\ast(0)}_2(z;\nu),
\end{equation}
\begin{equation}\label{eq:15.4.hr2f3}
\beta^{\ast(r+1)}_1(z;\nu)=\delta\beta^{\ast(r)}_2(z;\nu)-
\beta^{\ast(r)}_2(z;\nu)=
\end{equation}
$$\delta^{r+1}\beta^{(0)}_1(z;\nu)-\delta^r\beta^{(0)}_2(z;\nu),$$
where $r=0,\,1,\,2,$
and we obtain (\ref{eq:15.4.hr3}) -- (\ref{eq:15.4.hr3}) again. Moreover,
\begin{equation}\label{eq:15.4.hr2f4}
\beta^{\ast(r+1)}_4(z;\nu)=
\delta\beta^{\ast(r)}_4(z;\nu)+
\beta^{\ast(r)}_3(z;\nu),
\end{equation}
where $r=0,\,1,\,2.$
If we take now $z\in(1,+\infty)$ and will tend $z$ to $1,$ then,
in view of (\ref{eq:15.4.hr1f}), (\ref{eq:6cd1}), (\ref{eq:6cd2})
 and (\ref{eq:6ch0})
\begin{equation}\label{eq:15.4.hr6}
\delta^rf_{\alpha,0,2+j}^{\ast}(1,\nu)
=\lim\limits_{z\to1+0}\delta^rf_2^{\ast}(z, \nu)=
\end{equation}
$$(1-j+ij)\beta^{\ast(r)}_2(1;\nu)\zeta(2+j)-
\beta^{\ast(r)}_{3+j}(1;\nu)=$$
$$(1+j)\beta^{\ast(r)}_2(1;\nu)\zeta(2+j)-
\beta^{\ast(r)}_{3+j}(1;\nu),$$
where $r=0,\,1,\,2,\,i=2,\,j=0,\,1,$
\begin{equation}\label{eq:15.4.hr7}
\lim\limits_{z\to1+0} (z-1)\delta^3f_2^{\ast}(z, \nu)=0.
\end{equation}
In view of (\ref{eq:6cd}), (\ref{eq:6ce}),
 (\ref{eq:r0a1i2k01n0} -- (\ref{eq:r0a1i1k2n1}),
\begin{equation}\label{eq:a1r0i1n0}
\beta^{\ast(0)}_1(z;0)=\beta_{1,0,0}^{(0)}+
\beta_{1,1,0}^{(0)}z=-2+2z,
\end{equation}
\begin{equation}\label{eq:a1r0i2n0}
\beta^{\ast(0)}_2(z;0)=\beta_{2,0,0}^{(0)}+
\beta_{2,1,0}^{(0)}z=1+z,
\end{equation}
\begin{equation}\label{eq:a1r0i3n0}
\beta^{\ast(0)}_3(z;0)=\beta_{1,1,0}^{(0)}+
\beta_{2,1,0}^{(0)}=3,
\end{equation}
\begin{equation}\label{eq:a1r0i4n0}
\beta^{\ast(0)}_4(z;0)=\beta_{1,1,0}^{(0)}+
2\beta_{2,1,0}^{(0)}=4,
\end{equation}
\begin{equation}\label{eq:a1r0i1n1}
\beta^{\ast(0)}_1(z;1)=\beta_{1,0,1}^{(0)}+\beta_{1,1,1}^{(0)}z+
\end{equation}
$$\beta_{1,2,1}^{(0)}z^2=-5-16z+21z^2=(z-1)(21z+5),$$
\begin{equation}\label{eq:a1r0i2n1}
\beta^{\ast(0)}_2(z;1)=\beta_{2,0,1}^{(0)}+
\beta_{2,1,1}^{(0)}z+\beta_{2,2,1}^{(0)}z^2=1+16z+9z^2,
\end{equation}
\begin{equation}\label{eq:a1r0i3n1}
\beta^{\ast(0)}_3(z;1)=\beta_{1,1,1}^{(0)}+\beta_{2,1,1}^{(0)}+
\end{equation}
$$\frac12\beta_{1,2,1}^{(0)}+\frac14\beta_{2,2,1}^{(0)}+
(\beta_{1,2,1}^{(0)}+\beta_{2,2,1}^{(0)})z=$$
$$-16+16+\frac12\times21+\frac14\times9+(21+9)z=\frac{51}4+30z,$$
\begin{equation}\label{eq:a1r0i4n1}
\beta^{\ast(0)}_4(z;1)=\beta_{1,1,1}^{(0)}+2\beta_{2,1,1}^{(0)}+
\end{equation}
$$\frac14\beta_{1,2,1}^{(0)}+2\times\frac18\beta_{2,2,1}^{(0)}+
(\beta_{1,2,1}^{(0)}+2\beta_{2,2,1}^{(0)})z=$$
$$-16+2\times16+\frac14\times21+\frac14\times9+(21+18)z=\frac{47}2+39z.$$
\begin{equation}\label{eq:a1r0i1n2}
\beta^{\ast(0)}_1(z;2)=\beta_{1,0,2}^{(0)}+\beta_{1,1,2}^{(0)}z+
\end{equation}
$$\beta_{1,2,2}^{(0)}z^2+\beta_{1,3,2}^{(0)}z^3=$$
$$-\frac{20}3-216z-54z^2+\frac{830}3z^3=(z-1)(830z^2+668z+20)/3.$$
\begin{equation}\label{eq:a1r0i2n2}
\beta^{\ast(0)}_2(z;2)=\beta_{2,0,2}^{(0)}+\beta_{2,1,2}^{(0)}z+
\beta_{2,2,2}^{(0)}z^2+
\end{equation}
$$\beta_{2,3,2}^{(0)}z^3=1+81z+324z^2+100z^3.$$
\begin{equation}\label{eq:a1r0i3n2}
\beta^{\ast(0)}_{3}(z;2)=\beta_{1,1,2}^{(0)}+\beta_{2,1,2}^{(0)}+
\end{equation}
$$\frac12\beta_{1,2,2}^{(0)}+\frac14\beta_{2,2,2}^{(0)}+
\frac13\beta_{1,3,2}^{(0)}+\frac19\beta_{2,3,2}^{(0)}+$$
$$(\beta_{1,2,2}^{(0)}+\beta_{2,2,2}^{(0)})z+
\left(\frac12\beta_{1,3,2}^{(0)}+\frac14\beta_{2,3,2}^{(0)}\right)z+
(\beta_{1,3,2}^{(0)}+\beta_{2,3,2}^{(0)})z^2=
$$
%$$$$
$$-216+81-\frac{54}2+\frac{324}4+
\frac{830+100}9+(324-54+415/3+100/4)z+$$
$$(830/3+100)z^2=\frac{67}3+\frac{1300}3z+\frac{1130}3z^2.$$
\begin{equation}\label{eq:a1r0i4n2}
\beta^{\ast(0)}_4(z;2)=\beta_{1,1,2}^{(0)}+2\beta_{2,1,2}^{(0)}+
(1/4)\beta_{1,2,2}^{(0)}+2\times(1/8)\beta_{2,2,2}^{(0)}+
\end{equation}
%$$$$
$$(1/9)\beta_{1,3,2}^{(0)}+2\times(1/27)\beta_{2,3,2}^{(0)}+
(\beta_{1,2,2}^{(0)}+2\beta_{2,2,2}^{(0)})z+
$$
$$
\left(\frac14\beta_{1,3,2}^{(0)}+2\times\frac18\beta_{2,3,2}^{(0)}
\right)z+(\beta_{1,3,2}^{(0)}+2\beta_{2,3,2}^{(0)})z^2=-216+162+
$$
%$$$$
$$\frac{-54+324}4+\frac{830+200}{27}+((-57+648+(830/3+100)/4)z+$$
$$(830/3+200)z^2=2789/54+4129z/6+1430z^2/3.$$
In view of (\ref{eq:15.4.hr3}) -- (\ref{eq:15.4.hr5}), (\ref{eq:15.4.hr2f4})
and (\ref{eq:a1r0i1n0}) -- (\ref{eq:a1r0i4n2}),
\begin{equation}\label{eq:a1r1i1n0}
\beta^{\ast(1)}_1(z;0)=\delta\beta^{\ast(0)}_1(z;0)-
\beta^{\ast(0)}_2(z;0)=2z-1-z=z-1,
\end{equation}
\begin{equation}\label{eq:a1r1i2n0}
\beta^{\ast(1)}_2(z;0)=\delta\beta^{\ast(0)}_2(z;0)=z,
\end{equation}
\begin{equation}\label{eq:a1r1i3n0}
\beta^{\ast(1)}_3(z;0)=\delta\beta^{\ast(0)}_3(z;0)+
\beta^{\ast\vee(0)}_1(z;0)=2
\end{equation}
\begin{equation}\label{eq:a1r1i4n0}
\beta^{\ast(1)}_4(z;0)=\delta\beta^{\ast(0)}_4(z;0)+
\beta^{\ast(0)}_3(z;0)=3,
\end{equation}
\begin{equation}\label{eq:a1r1i1n1}
\beta^{\ast(1)}_1(z;1)=\delta\beta^{\ast(0)}_1(z;1)-
\beta^{\ast(0)}_2(z;1)=-16z+
\end{equation}
$$42z^2-(1+16z+9z^2)=-1-32z+33z^2=(z-1)(33z+1),$$
\begin{equation}\label{eq:a1r1i2n1}
\beta^{\ast(1)}_2(z;1)=\delta\beta^{\ast(0)}_2(z;1)=
16z+18z^2,
\end{equation}
\begin{equation}\label{eq:a1r1i3n1}
\beta^{\ast(1)}_3(z;1)=\delta\beta^{\ast(0)}_3(z;1)+
\beta^{\ast\vee(0)}_1(z;1)=30z+21z+5=51z+5,
\end{equation}
%$$$$
\begin{equation}\label{eq:a1r1i4n1}
\beta^{\ast(1)}_4(z;1)=\delta\beta^{\ast(0)}_4(z;1)+
\beta^{\ast(0)}_3(z;1)=39z+\frac{51}4+30z=69z+\frac{51}4,
\end{equation}
%$$$$
\begin{equation}\label{eq:a1r1i1n2}
\beta^{\ast(1)}_1(z;2)=\delta\beta^{\ast(0)}_1(z;2)-
\beta^{\ast(0)}_2(z;2)=
\end{equation}
$$-216z-108z^2+830z^3-(1+81z+324z^2+100z^3)=$$
$$-1-297z-432z^2+730z^3=(z-1)(730z^2+298z+1),$$
\begin{equation}\label{eq:a1r1i2n2}
\beta^{\ast(1)}_2(z;2)=\delta\beta^{\ast(0)}_2(z;2)=
81z+648z^2+300z^3,
\end{equation}
\begin{equation}\label{eq:a1r1i3n2}
\beta^{\ast(1)}_3(z;2)=\delta\beta^{\ast(0)}_3(z;2)+
\beta^{\ast\vee(0)}_1(z;2)=1300z/3+
\end{equation}
$$(2260z^2+830z^2+668z+20)/3=1030z^2+656z+\frac{20}3,$$
\begin{equation}\label{eq:a1r1i4n2}
\beta^{\ast(1)}_4(z;2)=\delta\beta^{\ast(0)}_4(z;2)+
\beta^{\ast(0)}_3(z;2)=4129z/6+2860z^2/3+
\end{equation}
$$(67+1300z+1130z^2)/3=67/3+6729z/6+1330z^2,$$
In view of (\ref{eq:15.4.hr3}) -- (\ref{eq:15.4.hr5}), (\ref{eq:15.4.hr2f4})
and (\ref{eq:a1r1i1n0}) -- (\ref{eq:a1r1i4n2}),
\begin{equation}\label{eq:a1r2i1n0}
\beta^{\ast(2)}_1(z;0)=\delta\beta^{\ast(1)}_1(z;0)-
\beta^{\ast(1)}_2(z;0)=z-z=0,
\end{equation}
\begin{equation}\label{eq:a1r2i2n0}
\beta^{\ast(2)}_2(z;0)=\delta\beta^{\ast(1)}_2(z;0)=z,
\end{equation}
\begin{equation}\label{eq:a1r2i3n0}
\beta^{\ast(2)}_3(z;0)=\delta\beta^{\ast(1)}_3(z;0)+
\beta^{\ast\vee(1)}_1(z;0)=1,
\end{equation}
\begin{equation}\label{eq:a1r2i4n0}
\beta^{\ast(2)}_4(z;0)=\delta\beta^{\ast(1)}_4(z;0)+
\beta^{\ast(1)}_3(z;0)=2
\end{equation}
\begin{equation}\label{eq:a1r2i1n1}
\beta^{\ast(2)}_1(z;1)=\delta\beta^{\ast(1)}_1(z;1)-
\beta^{\ast(1)}_2(z;1)=
\end{equation}
$$-32z+66z^2-(16z+18z^2)=-48z+48z^2=48z(z-1).$$
\begin{equation}\label{eq:a1r2i2n1}
\beta^{\ast(2)}_2(z;1)=\delta\beta^{\ast(1)}_2(z;1)=
16z+36z^2,
\end{equation}
\begin{equation}\label{eq:a1r2i3n1}
\beta^{\ast(2)}_3(z;1)=\delta\beta^{\ast(1)}_3(z;1)+
\beta^{\ast\vee(1)}_1(z;1)=
\end{equation}
$$51z+33z+1=84z+1,$$
\begin{equation}\label{eq:a1r2i4n1}
\beta^{\ast(2)}_4(z;1)=\delta\beta^{\ast(1)}_4(z;1)+
\beta^{\ast(1)}_3(z;1)=
\end{equation}
$$69z+51z+5=120z+5,$$
\begin{equation}\label{eq:a1r2i1n2}
\beta^{\ast(2)}_1(z;2)=\delta\beta^{\ast(1)}_1(z;2)-
\beta^{\ast(1)}_2(z;2)=
\end{equation}
$$-297z-864z^2+2190z^3-(81z+648z^2+300z^3)=$$
$$z(-378-1512z+1890z^2)=378(z-1)z(5z+1),$$
\begin{equation}\label{eq:a1r2i2n2}
\beta^{\ast(2)}_2(z;2)=\delta\beta^{\ast(1)}_2(z;2)=
81z+1296z^2+900z^3,
\end{equation}
\begin{equation}\label{eq:a1r2i3n2}
\beta^{\ast(0)}_3(z;2)=\delta\beta^{\ast(1)}_3(z;2)+
\beta^{\ast\vee(1)}_1(z;2)=
\end{equation}
$$2060z^2+656z+730z^2+298z+1=2790z^2+954z+1,$$
\begin{equation}\label{eq:a1r2i4n2}
\beta^{\ast(2)}_4(z;2)=\delta\beta^{\ast(1)}_4(z;2)+
\beta^{\ast(1)}_3(z;2)=6729z/6+2660z^2+
\end{equation}
$$1030z^2+656z+20/3=3690z^2+3555z/2+20/3,$$
%%%%%%%%%%%%%%%%%%%%%%%%%%%%%%%%%%%\S 5%%%%%%%%%%%%%%%%%%%%%%%%%%%%%%%%%%%%%%
\refstepcounter{section} {\begin{center}\large\bf \S 5. Auxiliary difference
equation.
 \end{center}}
Let $y_{i,k}(z;\nu)$ denotes $i-$th element of
 the column $Y_k(z;\nu)$
in (\ref{eq:15.1.h}). Then, in view of (\ref{eq:15e}), (\ref{eq:15.1.b}),
 (\ref{eq:15.1.h}),
\setcounter{equation}0
\begin{equation}\label{eq:7.1}
y_{j+1-\kappa,k}(z,\nu)=\delta^j f_k(z,\nu),\,y_{4,k}(z,\nu)=
\delta^3 f_k(z,\nu),
\end{equation}
where
$j=1,\,2,\,k=1,\,2,\,3,\,\vert z\vert>1,\,\nu\in{\mathbb N}_0.$
We denote $v^{\ast\ast}_{i,j}(\nu)$ the expression, which
stands in the matrix $V^{\ast\ast}(\nu)$ in intersection
of $i$-th row and $j$-th column, where $i=1,\,2,\,3,\,4,\,j=1,\,2,\,3,\,4.$
Let
\begin{equation}\label{eq:93bd1}
 D(z,\nu,w)=z(w^2-\mu)^{2}-w^{4},\,\mu=(\nu+1)^2
\end{equation}
In view of (\ref{eq:93bd})
\begin{equation}\label{eq:93bd1a}
(1/z)D(z,\nu,w)=(1-1/z)w^4+
\sum\limits_{k=0}^3r_{k+1}(\nu)w^k.
\end{equation}
It follows from general properties of Mejer's functions that
\begin{equation}\label{eq:93bd2}
 D(z,\nu,\delta)f_k(z,\nu)=0,
\end{equation}
where
$\vert z\vert>1,-3\pi/2<\arg(z)\le\pi/2,\log(z)=\ln(\vert z\vert)+i\arg(z),\,
k=1,\,2,\,3.$ Therefore,  in view of (\ref{eq:15e}), (\ref{eq:15.1.b}),
 (\ref{eq:15.1.h}),
\begin{equation}\label{eq:93bd3}
y_k(z,\nu)=
-(1-1/z)\delta^4 f_k(z,\nu)
\end{equation}
where
$$\vert z\vert>1,-3\pi/2<\arg(z)\le\pi/2,$$
$$\log(z)=\ln(\vert z\vert)+i\arg(z),\,k=1,\,2,\,3.$$
 In view of (\ref{eq:b}) -- (\ref{eq:ac}), (\ref{eq:acrf}),
\begin{equation}\label{eq:7.2}
\lim\limits_{z\to1+0}(z-1)\delta^4 f_2(z,\nu)=
\end{equation}
$$\lim\limits_{z\to1+0}(z-1)(O(1)\ln(1-1/z)+1/(z-1))=1,$$
\begin{equation}\label{eq:7.3}
\lim\limits_{z\to1+0}(z-1)\delta^4 f_k(z,\nu)=0,
\end{equation}
if $k-2=\pm1,$
\begin{equation}\label{eq:7.4}
\lim\limits_{z\to1+0}(\log(z))\delta^i f_k(z,\nu)=
\lim\limits_{z\to1+0}(z-1)\delta^i f_k(z,\nu)=0,
\end{equation}
%$$$$
if $i=0,\,1,\,2,\,3,k=1,\,2,\,3.$
Hence, if we tend $z\in(1,+\infty)$ to $1,$ then,
 in view of (\ref{eq:15e}), (\ref{eq:15.1.b}), we obtain
the equalities
\begin{equation}\label{eq:7.5}
y_{1,1}(1,\nu)=y_{1,3}(1,\nu)=0,y_{1,2}(1,\nu)=-1.
\end{equation}
In view of (\ref{eq:15.1.aj}), (\ref{eq:15.2.vw1,234}) -- (\ref{eq:a43asas}),
(\ref{eq:7.1}), (\ref{eq:93bd3}),
\begin{equation}\label{eq:7.6}
-a_{i,1}^{\ast\ast}(1;\nu)(1-1/z)\delta^4 f_k(z,\nu)+
\left(\sum\limits_{j=1}^2a_{i+1,j+1}^{\ast\ast}(1;\nu)
\delta^jf_{k}(z,\nu)\right)-\end{equation}
%$$$$
$$
(z-1)v_{i,1}^{\ast\ast}(\nu)(1-1/z)\delta^4 f_k(z,\nu)+
$$
$$(z-1)\sum\limits_{j=1}^2
v_{i+1,j+1}^{\ast\ast}(\nu)\delta^{j}f_{k}(z,\nu)=
%$$$$
\mu_1(\nu)^2\nu^5\delta^if_k(z,\nu-1),$$
where
$i=1,\,2,\,k=1,\,2,\,3,\,\vert z\vert>1,\,-3\pi/2<\arg(z)\le\pi/2$ and
$\nu$ run over the
 set $M_1^{\ast}=((-\infty,-2]\cup[1,+\infty))\cap{\mathbb Z}.$
We tend $z\in(1,+\infty)$ to $1$ now and obtain the equalities
\begin{equation}\label{eq:7.7}
a_{i+1,1}^{\ast\ast}(1;\nu)(k-1)(k-3)+
\left(\sum\limits_{j=1}^2a_{i+1,j+1}^{\ast\ast}(1;\nu)
(\delta^jf_k)(1,\nu)\right)=
\end{equation}
$$$$
$$\mu_1(\nu)^2\nu^5\delta^if_k(1,\nu-1),$$
where
$i=1,\,2,\,k=1,\,2,\,3$  and
$\nu\in M_1^{\ast}=((-\infty,-2]\cup[1,+\infty))\cap{\mathbb Z}.$
Let are given
\begin{equation}\label{eq:fg}
F=\{F(\nu)\}_{\nu=-\infty}^{+\infty}\,\,\text{and}\,\,
G=\{G(\nu)\}_{\nu=-\infty}^{+\infty}
\end{equation}
 such that
\begin{equation}\label{eq:fg1}
F(-\nu-2)=F(\nu),\,G(-\nu-2)=G(\nu)
,\,F(\nu)\in{\mathbb R},\,G(\nu)\in{\mathbb R}
\end{equation}
for $\nu\in{\mathbb Z}.$
 Let further
\begin{equation}\label{eq:yfg}
y_{F,G}^{\ast\ast}(z,\nu)=F(\nu)
\delta f_k(z,\nu)+
G(\nu)\delta^2f_k(z,\nu)
\end{equation}
for $k=1,\,2,\,3$ and
$\nu\in M_1^{\ast\ast\ast}=((-\infty,-2]\cup[0,+\infty))\cap{\mathbb Z}.$
Since $F$ and $G$ have the property (\ref{eq:fg1}),
it follows from (\ref{eq:15.1.i})) that
\begin{equation}\label{eq:yfg2}
y_{F,G}^{\ast\ast}(z,-\nu-2)=y_{F,G}^{\ast\ast}(z,\nu)
\end{equation}
for $k=1,\,2,\,3$ and
$\nu\in M_1^{\ast\ast\ast}=((-\infty,-2]\cup[0,+\infty))\cap{\mathbb Z}.$
Let
\begin{equation}\label{eq:afgkpj}
a_{F,G,j}^{\ast\ast\ast}(z;\nu)=F(\nu-1)a_{2,j}^{\ast\ast}(z;\nu)+
G(\nu-1)a_{3,j}^{\ast\ast}(z,\nu)
\end{equation}
for $\nu\in M_1^{\ast}=((-\infty,-2]\cup[1,+\infty))\cap{\mathbb Z}
,\,j=1,\,2,\,3.$ In view of (\ref{eq:7.7})
\begin{equation}\label{eq:7.7a}
a_{F,G,1}^{\ast\ast\ast}(1;\nu)(k-1)(k-3)+
\end{equation}
$$\left(\sum\limits_{j=1}^2
a_{F,G,j+1}^{\ast\ast\ast}(1;\nu)
(\delta^jf_{1,0,k})(1,\nu)\right)=%$$$$
\mu_1(\nu)^{2}\nu^5y_{F,G}^{\ast\ast}(1,\nu-1),$$
with $\nu\in M_1^{\ast}=((-\infty,-2]\cup[1,+\infty))\cap{\mathbb Z}
,\,k=1,\,2,\,3.$

Replacing  $\nu\in M_1^\ast$ by
$\nu:=-\nu-2\in M_1^{\ast\ast}=((-\infty,-3]\cup[0,+\infty))\cap{\mathbb Z}$

in (\ref{eq:7.7a}), and taking in account  (\ref{eq:15.1.i})
 we obtain the equalities
\begin{equation}\label{eq:7.8}
a_{F,G,1}^{\ast\ast\ast}(1;-\nu-2)(k-1)(k-3)+
\end{equation}
$$\left(\sum\limits_{j=1}^2
a_{F,G,j}^{\ast\ast\ast}(1;-\nu-2)
(\delta^jf_k)(1,\nu)\right)=-\mu_1(\nu)^{2}(\nu+2)^5
y_{F,G}^{\ast\ast}(z,\nu+1),$$
%$$$$
where $k=1,\,2,\,3$ and
$\nu\in M_1^{\ast\ast}=((-\infty,-3]\cup[0,+\infty))\cap{\mathbb Z}.$
Let
\begin{equation}\label{eq:7.12}
\vec w_{F,G,j}(\nu)=
\left(\begin{matrix} a_{F,G,j+1}^{\ast\ast\ast}(1;-\nu-2)\\
F(\nu)(2-j)+G(\nu)(j-1)\\
 a_{F,G,j+1}^{\ast\ast\ast}(1;\nu)
\end{matrix}\right),\end{equation}
where $j=1,2,\,
\nu\in M_1^{\ast\ast\ast\ast}=((-\infty,-3]\cup[1,+\infty))\cap{\mathbb Z},$
\begin{equation}\label{eq:7.13}
W_{F,G}(\nu)=\left(\begin{matrix} a_{F,G,2}^{\ast\ast\ast}(1;-\nu-2)&
a_{F,G,3}^{\ast\ast\ast}(1;-\nu-2)\\
F(\nu)&G(\nu)\\
 a_{F,G,2}^{\ast\ast\ast}(1;\nu)&
a_{F,G,3}^{\ast\ast\ast}(1;\nu)
\end{matrix}\right)=
\end{equation}
$$
\left(\begin{matrix} \vec w_{F,G,1}(\nu)&\vec w_{F,G,2}(\nu)\end{matrix}\right)
,\,Y^{\ast\ast\ast}_k(\nu)=
%$$$$
\left(\begin{matrix} (\delta f_k)(1,\nu)\\
(\delta^2f_k)(1,\nu)\end{matrix}\right),$$
\begin{equation}\label{eq:7.14}
Y^{\ast\ast\ast\ast}_{F,G,k}(\nu)=
\left(\begin{matrix}
\mu_1(-\nu-2)^2(-\nu-2)^5y_{F,G}^{\ast\ast}(z,-\nu-3)\\
y_{F,G}^{\ast\ast}(z,\nu)\\
\mu_1(\nu)^2\nu^5y_{F,G}^{\ast\ast}(z,\nu-1)
\end{matrix}\right),
\end{equation}
%$$$$
where $k=1,\,3,\,
\nu\in M_1^{\ast\ast\ast\ast}=((-\infty,-3]\cup[1,+\infty))\cap{\mathbb Z}.$
Let further
\begin{equation}\label{eq:7.11}
\vec w_{F,G,3}(\nu)=
\left(\begin{matrix} w_{F,G,3,1}(\nu)\\
                w_{F,G,3,2}(\nu)\\
                 w_{F,G,3,3}(\nu)
\end{matrix}\right)=[\vec w_{F,G,1}(\nu),\vec w_{F,G,2}(\nu)].
\end{equation}
is vector product of $\vec w_{F,G,1}(\nu)$ and $\vec w_{F,G,2}(\nu).$

Let $\bar w_{F,G,3}(\nu)=(\vec w_{F,G,3}(\nu))^t$ is the row conjugate to the
column $\vec w_{F,G,3}(\nu).$ Then for scalar products $(\vec
w_{F,G,3}(\nu),\vec w_{F,G,j}(\nu))$ we have the equalities
$$\bar w_{F,G,3}(\nu)\vec w_{F,G,j}(\nu)=
(\vec w_{F,G,3}(\nu),\vec w_{F,G,j}(\nu))=0,$$
where
$\nu\in M_1^{\ast\ast\ast\ast}=((-\infty,-3]\cup[1,+\infty))\cap{\mathbb Z},\,
j=1,\,2.$ Therefore
\begin{equation}\label{eq:7.9}
\bar w_{F,G,3}(\nu)W_{F,G}(\nu)=
\left(\begin{matrix} 0&0\end{matrix}\right),
\end{equation}
where
$\nu\in M_1^{\ast\ast\ast\ast}=((-\infty,-3]\cup[1,+\infty))\cap{\mathbb Z}.$

In view of (\ref{eq:7.7}) (\ref{eq:7.8}) and (\ref{eq:7.9}),
\begin{equation}\label{eq:7.10}
\bar w_{F,G,3}(\nu)Y^{\ast\ast\ast\ast}_{F,G,k}(\nu)=
\bar w_{F,G,3}(\nu)W_{F,G,3}(\nu)Y^{\ast\ast\ast}_k(\nu)=0,
\end{equation}
%$$$$
where $k=1,\,3$ and
$\nu\in M_1^{\ast\ast\ast\ast}=((-\infty,-3]\cup[1,+\infty))\cap{\mathbb Z}.$

In view of (\ref{eq:afgkpj}), (\ref{eq:7.11}) -- (\ref{eq:7.13}),
% (\ref{eq:15.2.e}),
 (\ref{eq:15.2.vw1,234}) -- (\ref{eq:a43asas}), and,

 since $\tau_{-\nu-2}=-\nu-1=-\tau_\nu=-\tau,$ it follows that
\begin{equation}\label{eq:afg01t}
a_{F,G,1}^{\ast\ast\ast}(1;\nu)=F(\nu-1)a_{2,1}^{\ast\ast}(1,\nu)+
G(\nu-1)a_{3,1}^{\ast\ast}(1,\nu)=
\end{equation}
%$$$$
$$-F(\nu-1)\tau^2(\tau-1)(2\tau-1)(6\tau^2-4\tau+1)+$$
$$G(\nu-1)\tau^2(\tau-1)^2(2\tau-1)(4\tau^2-3\tau+1),$$
\begin{equation}\label{eq:afg01mt}
a_{F,G,1}^{\ast\ast\ast}(1;-\nu-2)=
\end{equation}
$$F(\nu+1)a_{2,1}^{\ast\ast}(1,-\nu-2)+
G(\nu+1)a_{3,1}^{\ast\ast}(1,-\nu-2)=$$
$$-F(\nu+1)\tau^2(\tau+1)(2\tau+1)(6\tau^2+4\tau+1)-$$
$$-G(\nu+1)\tau^2(\tau+1)^2(2\tau+1)(4\tau^2+3\tau+1),$$
\begin{equation}\label{eq:afg02t}
a_{F,G,2}^{\ast\ast\ast}(1;\nu)=
\end{equation}
$$F(\nu-1)a_{1,0,2,2}^{\ast\ast}(1,\nu)+
G(\nu-1)a_{3,2}^{\ast\ast}(1,\nu)=$$
$$F(\nu-1)\tau^5(\tau-1)(\tau^3+2(2\tau-1)^3)-$$
$$2G(\nu-1)\tau^5(\tau-1)^2(2\tau-1)(\tau^3-(\tau-1)^3),$$
\begin{equation}\label{eq:afg02mt}
a_{F,G,2}^{\ast\ast\ast}(z;-\nu-2)=
\end{equation}
$$F(\nu+1)a_{2,2}^{\ast\ast}(1,-\nu-2)+
G(\nu+1)a_{3,2}^{\ast\ast}(1,-\nu-2)=$$
$$-F(\nu+1)\tau^5(\tau+1)(\tau^3+2(2\tau+1)^3)-$$
$$G(\nu+1)2\tau^5(\tau+1)^2(2\tau+1)((\tau+1)^3-\tau^3),$$
\begin{equation}\label{eq:afg03t}
a_{F,G,3}^{\ast\ast\ast}(1;\nu)=
\end{equation}
$$F(\nu-1)a_{2,3}^{\ast\ast}(1,\nu)+
G(\nu-1)a_{3,3}^{\ast\ast}(1,\nu)=$$
$$-3F(\nu-1)\tau^4(\tau-1)(2\tau-1)^3+
G(\nu-1)\tau^4(\tau-1)^2((\tau-1)^3+2(2\tau-1)^3),$$
\begin{equation}\label{eq:afg03mt}
a_{F,G,3}^{\ast\ast\ast}(z;-\nu-2)=
\end{equation}
$$F(\nu+1)a_{2,3}^{\ast\ast}(1,-\nu-2)+
G(\nu+1)a_{3,3}^{\ast\ast}(1,-\nu-2)=$$
$$-3F(\nu+1)\tau^4(\tau+1)(2\tau+1)^3-$$
$$G(\nu+1)\tau^4(\tau+1)^2(2\tau+1)((\tau+1)^3+2(2\tau+1)^3),$$
%%%%%%%%%%%%%%%%%%%%%%%%%%%%%F,G-constants%%%%%%%%%%%%%%%%%%%%%%%%%%%%%
We consider the case now, when $F$ and $G$ are constant sequences, or,
equivalently, real constants.
%such that
%\begin{equation}\label{eq:FG}
% F+G\ne0,\,(F+G)G\ge0.
%\end{equation}
In view of (\ref{eq:7.12}), (\ref{eq:afg02t}) -- (\ref{eq:afg03mt})
\begin{equation}\label{eq:7.12F1G0j}
\vec w_{1,0,j}(\nu)=
\left(\begin{matrix} a_{1,0,j+1}^{\ast\ast\ast}(1;-\nu-2)\\
2-j\\
 a_{1,0,j+1}^{\ast\ast\ast}(1;\nu),
\end{matrix}\right),\end{equation}
where $j=1,2,\,
\nu\in M_1^{\ast\ast\ast\ast}=((-\infty,-3]\cup[1,+\infty))\cap{\mathbb Z},$
\begin{equation}\label{eq:7.12F1G0j1}
\vec w_{1,0,1}(\nu)=
\left(\begin{matrix} a_{1,0,2}^{\ast\ast\ast}(1;-\nu-2)\\
1\\
 a_{1,0,2}^{\ast\ast\ast}(1;\nu)
\end{matrix}\right)=
\end{equation}
$$\left(\begin{matrix}
-\tau^5(\tau+1)(\tau^3+2(2\tau+1)^3)\\
1\\
\tau^5(\tau-1)(\tau^3+2(2\tau-1)^3)
\end{matrix}\right),$$
\begin{equation}\label{eq:7.12F1G0j2}
\vec w_{1,0,2}(\nu)=
\left(\begin{matrix} a_{1,0,3}^{\ast\ast\ast}(1;-\nu-2)\\
0\\
 a_{1,0,3}^{\ast\ast\ast}(1;\nu)
\end{matrix}\right)=\end{equation}
$$\left(\begin{matrix}
-3\tau^4(\tau+1)(2\tau+1)^3\\
0\\
-3\tau^4(\tau-1)(2\tau-1)^3
\end{matrix}\right),$$
\begin{equation}\label{eq:7.12F0G1j}
\vec w_{0,1,j}(\nu)=
\left(\begin{matrix} a_{0,1,j+1}^{\ast\ast\ast}(1;-\nu-2)\\
j-1\\
 a_{0,1,j+1}^{\ast\ast\ast}(1;\nu),
\end{matrix}\right),\end{equation}
where $j=1,2,\,
\nu\in M_1^{\ast\ast\ast\ast}=((-\infty,-3]\cup[1,+\infty))\cap{\mathbb Z},$
\begin{equation}\label{eq:7.12F0G1j1}
\vec w_{0,1,1}(\nu)=
\left(\begin{matrix} a_{0,1,2}^{\ast\ast\ast}(1;-\nu-2)\\
0\\
 a_{0,1,2}^{\ast\ast\ast}(1;\nu)
\end{matrix}\right)=
\end{equation}
$$\left(\begin{matrix}
-2\tau^5(\tau+1)^2(2\tau+1)((\tau+1)^3-\tau^3) \\
0\\
-2\tau^5(\tau-1)^2(2\tau-1)(\tau^3-(\tau-1)^3)
\end{matrix}\right),$$
\begin{equation}\label{eq:7.12F0G1j2}
\vec w_{0,1,2}(\nu)=
\left(\begin{matrix} a_{0,1,3}^{\ast\ast\ast}(1;-\nu-2)\\
1\\
 a_{0,1,3}^{\ast\ast\ast}(1;\nu)
\end{matrix}\right)=
\end{equation}
$$\left(\begin{matrix}
-\tau^4(\tau+1)^2(2\tau+1)((\tau+1)^3+2(2\tau+1)^3)\\
1\\
\tau^4(\tau-1)^2((\tau-1)^3+2(2\tau-1)^3)
\end{matrix}\right),$$
and,since $F,\,G$ are constants now, it follows that
\begin{equation}\label{eq:7.12h}
\vec w_{F,G,j}(\nu)=F\vec w_{1,0,j}(\nu)+G\vec w_{0,1,j}(\nu)
\end{equation}
where $j=1,2,\,
\nu\in M_1^{\ast\ast\ast\ast}=((-\infty,-3]\cup[1,+\infty))\cap{\mathbb Z}.$
In view of (\ref{eq:7.11}),
\begin{equation}\label{eq:7.11h}
\vec w_{F,G,3}(\nu)=F^2[\vec w_{1,0,1}(\nu),\vec w_{1,0,2}(\nu)]+
\end{equation}
$$FG([\vec w_{1,0,1}(\nu),\vec w_{0,1,2}(\nu)]+
[\vec w_{0,1,1}(\nu),\vec w_{1,0,2}(\nu)])+%$$$$
G^2[\vec w_{0,1,1}(\nu),\vec w_{0,1,2}(\nu)].$$
{\bf For any
$$\vec a=\left(\begin{matrix} a_1\\a_2\\a_3\end{matrix}\right)$$
we put $(\vec a)_i=a_i$ for $i=1,\,2,\,3.$}
Let further
\begin{equation}\label{eq:7.11hh}
\vec w_{i,j,4}(\nu)=
%\left(\matrix w_{i,j,4,1}(\nu)\\
%                w_{i,j,4,2}(\nu)\\
%                 w_{i,j,4,3}(\nu)
%\endmatrix\right)=
[\vec w_{i,1-i,1}(\nu),\vec w_{j,1-j,2}(\nu)],
\end{equation}
with $i=0,\,1,\,j=0,\,1.$
%%%%%%%%%%%%%%%%%%%%%%%%%%%%%%%w}_{1,1,4,1}(\nu)%%%%%%%%%%%%%%%%%%%%%%
In view of (\ref{eq:7.11hh}),(\ref{eq:afg02t})-- (\ref{eq:afg03mt}),
$$\vec w_{1,1,4}(\nu)=
[\vec w_{1,0,1}(\nu),\vec w_{1,0,2}(\nu)]=\vec w_{1,0,3}(\nu)],$$
\begin{equation}\label{eq:7w01031}
(\vec w_{1,1,4}(\nu))_1=(\vec w_{1,0,3}(\nu))_1=
\end{equation}
$$\det\left(\begin{matrix}
1&0\\
 a_{1,0,2}^{\ast\ast\ast}(1;\nu)&
a_{1,0,3}^{\ast\ast\ast}(1;\nu)
\end{matrix}\right)=
$$
$$a_{1,0,3}^{\ast\ast\ast}(1;\nu)=a_{2,3}^{\ast\ast}(1,\nu)
=-3\tau^4(\tau-1)(2\tau-1)^3,$$
%%%%%%%%%%%%%%%%%%%%%%%%%%%%%%%%w_{1,1,4,2}(\nu)%%%%%%%%%%%%%%%%%%%
\begin{equation}\label{eq:7w01032}
(\vec w_{1,1,4}(\nu))_2=(\vec w_{1,0,3}(\nu))_2=
\end{equation}
$$-\det\left(\begin{matrix}
 a_{1,0,2}^{\ast\ast\ast}(1;-\nu-2)&
a_{1,0,3}^{\ast\ast\ast}(1;-\nu-2)
\\
 a_{1,0,2}^{\ast\ast\ast}(1;\nu)&
a_{1,0,3}^{\ast\ast\ast}(1;\nu)
\end{matrix}\right)=$$
$$-\det\left(\begin{matrix}
 a_{2,2}^{\ast\ast}(1;-\nu-2)&
a_{2,3}^{\ast\ast}(1;-\nu-2)
\\
 a_{2,2}^{\ast\ast}(1;\nu)&
a_{2,3}^{\ast\ast}(1;\nu)
\end{matrix}\right)=$$
$$a_{2,2}^{\ast\ast}(1;\nu)a_{2,3}^{\ast\ast}(1;-\nu-2)-
a_{2,3}^{\ast\ast}(1;\nu)a_{2,3}^{\ast\ast}(1;-\nu-2)=$$
$$\tau^5(\tau-1)(\tau^3+2(2\tau-1)^3)
(-3\tau^4(\tau+1)(2\tau+1)^3)-$$
$$(-3\tau^4(\tau-1)(2\tau-1)^3)
(-\tau^5(\tau+1)(\tau^3+2(2\tau+1)^3))=$$
$$-3\tau^9(\tau^2-1)
(\tau^3((2\tau-1)^3+(2\tau+1)^3)+4(4\tau^2-1)^3)=$$
$$-12\tau^9(\tau^2-1)(68\tau^6-45\tau^4+12\tau^2-1),$$
%%%%%%%%%%%%%%%%%%%%%%%%%%%%%%w_{1,1,4,3}(\nu)%%%%%%%%%%%%%%%%%%%%
\begin{equation}\label{eq:7w01033}
(\vec w_{1,1,4}(\nu))_3=(\vec w_{1,0,3}(\nu))_3=
\end{equation}
$$\det\left(\begin{matrix}
 a_{1,0,2}^{\ast\ast\ast}(1;-\nu-2)&
a_{1,0,3}^{\ast\ast\ast}(1;-\nu-2)\\
 1&0
\end{matrix}\right)=$$
$$
-a_{1,0,3}^{\ast\ast\ast}(1;-\nu-2)=-a_{2,3}^{\ast\ast}(1,-\nu-2)=
=3\tau^4(\tau+1)(2\tau+1)^3.
$$
%%%%%%%%%%%%%%%%%%%%%%%%w_{0,0,4,1}(\nu)%%%%%%%%%%%%%%%%%%%%%%%%%%%%
In view of (\ref{eq:7.11hh}),(\ref{eq:afg02t})-- (\ref{eq:afg03mt}),
$$\vec w_{0,0,4}(\nu)=
[\vec w_{0,1,1}(\nu),\vec w_{0,1,2}(\nu)]=\vec w_{0,1,3}(\nu),$$
\begin{equation}\label{eq:7w00131}
(\vec w_{0,0,4}(\nu))_1=(\vec w_{0,1,3}(\nu))_1=
\end{equation}
$$\det\left(\begin{matrix}
0&1\\
 a_{0,1,2}^{\ast\ast\ast}(1;\nu)&
a_{0,1,3}^{\ast\ast\ast}(1;\nu)
\end{matrix}\right)=$$
$$-a_{0,1,2}^{\ast\ast\ast}(1;\nu)=-a_{3,2}^{\ast\ast}(1,\nu)=
2\tau^5(\tau-1)^2(2\tau-1)(\tau^3-(\tau-1)^3),$$
%%%%%%%%%%%%%%%%%%%%%%%%%%%%%%%%w_{0,0,4,2}(\nu)%%%%%%%%%%%%%%%%%%%
\begin{equation}\label{eq:7w00132}
(\vec w_{0,0,4}(\nu))_2=(\vec w_{0,1,3}(\nu))_2=
\end{equation}
$$-\det\left(\begin{matrix}
 a_{0,1,2}^{\ast\ast\ast}(1;-\nu-2)&
a_{0,1,3}^{\ast\ast\ast}(1;-\nu-2)
\\
 a_{0,1,2}^{\ast\ast\ast}(1;\nu)&
a_{0,1,3}^{\ast\ast1\ast}(1;\nu)
\end{matrix}\right)=$$
$$-\det\left(\begin{matrix}
 a_{3,2}^{\ast\ast}(1;-\nu-2)&
a_{3,3}^{\ast\ast}(1;-\nu-2)
\\
 a_{3,2}^{\ast\ast}(1;\nu)&
a_{3,3}^{\ast\ast}(1;\nu)
\end{matrix}\right)=$$
$$a_{3,2}^{\ast\ast}(1;\nu)a_{3,3}^{\ast\ast}(1;-\nu-2)-
a_{3,3}^{\ast\ast}(1;\nu)a_{3,2}^{\ast\ast}(1;-\nu-2)=$$
$$-2\tau^5(\tau-1)^2(2\tau-1)(\tau^3-(\tau-1)^3)\times$$
$$(-\tau^4(\tau+1)^2((\tau+1)^3+2(2\tau+1)^3)-$$
$$(-2\tau^5(\tau+1)^2(2\tau+1)((\tau+1)^3-\tau^3))\times$$
$$(\tau^4(\tau-1)^2((\tau-1)^3+2(2\tau-1)^3)=$$
$$4\tau^9(\tau^2-1)^2(102\tau^6-68\tau^4+21\tau^2-3),$$
%%%%%%%%%%%%%%%%%%%%%%%%%%%%%%%%%%%%%%%%%%%%%%%%%%%%%%%%%%%%%%%%%%%%%%%%%%%%
%%%%%%%%%%%%%%%%%%%%%%%%%%%%%%%%w_{0,0,4,3}(\nu)%%%%%%%%%%%%%%%%%%%
\begin{equation}\label{eq:7w00133}
(\vec w_{0,0,4}(\nu))_3=(\vec w_{0,1,3}(\nu))_3=
\end{equation}
$$\det\left(\begin{matrix}
 a_{0,1,2}^{\ast\ast\ast}(1;-\nu-2)&
a_{0,1,3}^{\ast\ast\ast}(1;-\nu-2)\\
0&1
\end{matrix}\right)=$$
$$a_(3,2)^{\ast\ast}(1;-\nu-2)=
-2\tau^5(\tau+1)^2(2\tau+1)((\tau+1)^3-\tau)^3),$$
%%%%%%%%%%%%%%%%%%%%%%%%%%%%%%%%%%%%%%%%%%%%%%%%%%%%%%%%%%%%%%%%%%%%%%%%%%%%%%
%%%%%%%%%%%%%%%%%%%%%%%%w_{0,1,4,1}(\nu)%%%%%%%%%%%%%%%%%%%%%%%%%%%%
In view of (\ref{eq:7.11hh}),(\ref{eq:afg02t})-- (\ref{eq:afg03mt}),
$$\vec w_{0,1,4}(\nu)=
[\vec w_{0,1,1}(\nu),\vec w_{1,0,2}(\nu)],$$
\begin{equation}\label{eq:7w00141}
(\vec w_{0,1,4}(\nu))_1=([\vec w_{0,1,1}(\nu),\vec w_{1,0,2}(\nu)])_1=
\end{equation}
$$\det\left(\begin{matrix}
0&0\\
 a_{0,1,2}^{\ast\ast\ast}(1;\nu)&
a_{1,0,3}^{\ast\ast\ast}(1;\nu)
\end{matrix}\right)=0,$$
%%%%%%%%%%%%%%%%%%%%%%%%%%%%%%%%w_{0,1,4,2}(\nu)%%%%%%%%%%%%%%%%%%%
\begin{equation}\label{eq:7w00142}
(\vec w_{0,1,4}(\nu))_2=([\vec w_{0,1,1}(\nu),\vec w_{1,0,2}(\nu)])_2=
\end{equation}
$$-\det\left(\begin{matrix}
 a_{0,1,2}^{\ast\ast\ast}(1;-\nu-2)&
a_{1,0,3}^{\ast\ast\ast}(1;-\nu-2)
\\
 a_{0,1,2}^{\ast\ast\ast}(1;\nu)&
a_{1,0,3}^{\ast\ast\ast}(1;\nu)
\end{matrix}\right)=$$
$$-a_{3,2}^{\ast\ast}(1;-\nu-2)a_{2,3}^{\ast\ast}(1;\nu)+
a_{3,2}^{\ast\ast}(1;\nu)a_{2,3}^{\ast\ast}(1;-\nu-2)=$$
$$-12t^9(\tau^2-1)(4\tau^2-1)(12\tau^4-6\tau^2+1),$$
%%%%%%%%%%%%%%%%%%%%%%%%%%%%%%%%%%%%%%%%%%%%%%%%%%%%%%%%%%%%%%%%%%%%%%%%%%%%%%
%%%%%%%%%%%%%%%%%%%%%%%%%%%%%%%%w_{0,1,4,3}(\nu)%%%%%%%%%%%%%%%%%%%
\begin{equation}\label{eq:7w00143}
(\vec w_{0,1,4}(\nu))_3=([\vec w_{0,1,1}(\nu),\vec w_{1,0,2}(\nu)])_3,
\end{equation}
$$\det\left(\begin{matrix}
 a_{0,1,2}^{\ast\ast\ast}(1;-\nu-2)&
a_{1,0,3}^{\ast\ast\ast}(1;-\nu-2)
\\
0&0
\end{matrix}\right)=0,$$
%%%%%%%%%%%%%%%%%%%%%%%%%%%%%%%%%%%%%%%%%%%%%%%%%%%%%%%%%%%%%%%%%%%%%%%%%%%%%%
%%%%%%%%%%%%%%%%%%%%%%%%w_{1,0,4,1}(\nu)%%%%%%%%%%%%%%%%%%%%%%%%%%%%
In view of (\ref{eq:7.11hh}),(\ref{eq:afg02t})-- (\ref{eq:afg03mt}),
$$\vec w_{1,0,4}(\nu)=
[\vec w_{1,0,1}(\nu),\vec w_{0,1,2}(\nu)],$$
\begin{equation}\label{eq:7w01041}
(\vec w_{1,0,4}(\nu))_1=([\vec w_{1,0,1}(\nu),\vec w_{0,1,2}(\nu)])_1
\end{equation}
$$\det\left(\begin{matrix}
1&1\\
 a_{1,0,2}^{\ast\ast\ast}(1;\nu)&
a_{0,1,3}^{\ast\ast\ast}(1;\nu)
\end{matrix}\right)=$$
$$ a_{3,3}^{\ast\ast}(1;\nu)- a_{2,2}^{\ast\ast}(1;\nu)=$$
$$-t^4(t-1)(2t-1)(10t^2-10t+3),$$
%%%%%%%%%%%%%%%%%%%%%%%%w_{1,0,4,2}(\nu)%%%%%%%%%%%%%%%%%%%%%%%%%%%%
\begin{equation}\label{eq:7w01042}
(\vec w_{1,0,4}(\nu))_2=([w_{1,0,1}(\nu),w_{0,1,2}(\nu)])_2=
\end{equation}
$$-\det\left(\begin{matrix}
 a_{1,0,2}^{\ast\ast\ast}(1;-\nu-2)&
a_{0,1,3}^{\ast\ast\ast}(1;-\nu-2)\\
 a_{1,0,2}^{\ast\ast\ast}(1;\nu)&
a_{0,1,3}^{\ast\ast\ast}(1;\nu)
\end{matrix}\right)=$$
$$-a_{2,2}^{\ast\ast}(1;-\nu-2)a_{3,3}^{\ast\ast}(1;\nu)
+a_{2,2}^{\ast\ast}(1;\nu)a_{3,3}^{\ast\ast}(1;-\nu-2)=$$
$$-4t^9(t^2-1)(170t^6-104t^4+30t^2-3),$$
%%%%%%%%%%%%%%%%%%%%%%%%w_{1,0,4,3}(\nu)%%%%%%%%%%%%%%%%%%%%%%%%%%%%
\begin{equation}\label{eq:7w01043}
(\vec w_{1,0,4}(\nu))_3=([\vec w_{1,0,1}(\nu),\vec w_{0,1,2}(\nu)])_3=
\end{equation}
$$\det\left(\begin{matrix}
 a_{1,0,2}^{\ast\ast\ast}(1;-\nu-2)&
a_{0,1,3}^{\ast\ast\ast}(1;-\nu-2)\\
1&1
\end{matrix}\right)=$$
$$a_{2,2}^{\ast\ast}(1;-\nu-2)
-a_{3,3}^{\ast\ast}(1;-\nu-2)=$$
$$t^4(t+1)(2t+1)(10t^2+10t+3),$$
%%%%%%%%%%%%%%%%%%%%%%%%%%%%%%%%%%%%%%%%%%%%%%%%%%%%%%%%%%%%%%%%%%%%%%%%
%%%%%%%%%%%%%w_{1,0,4,2}(\nu)+w_{0,1,4,2}(\nu)%%%%%%%%%%%%%%%%%%%%%%%%%%%%
\begin{equation}\label{eq:7w010p0142}
(\vec w_{0,1,4}(\nu))_2+(\vec w_{1,0,4}(\nu))_2=
\end{equation}
$$-8t^9(t^2-1)(157t^6-106t^4+30t^2-3).$$
Therefore,
\begin{equation}\label{eq:7wkpfg31c}
(\vec w_{F,G,3}(\nu))_1=-3\tau^4(\tau-1)(2\tau-1)^3F^2-
\end{equation}
$$t^4(\tau-1)(2\tau-1)(10\tau^2-10\tau+3)FG+$$
$$2\tau^5(\tau-1)^2(2\tau-1)(\tau^3-(\tau-1)^3)G^2=$$
$$\tau^4(\tau-1)(2\tau-1)\times$$
$$(-3(2\tau-1)^2F^2-(10\tau^2-10\tau+3)FG+
2(\tau-1)(3\tau^3-3\tau^2+\tau)G^2)=$$
$$\tau^4(\tau-1)(2\tau-1)\times$$
$$(-3(2\tau-1)^2F^2-(10\tau^2-10\tau+3)FG+2(3\tau^4-6\tau^3+4\tau^2-\tau)G^2,$$
\begin{equation}\label{eq:7wkpfg32c}
(\vec w_{F,G,3}(\nu))_2=-12\tau^9(\tau^2-1)(68\tau^6-45\tau^4+12\tau^2-1)F^2-
\end{equation}
$$8t^9(t^2-1)(157t^6-106t^4+30t^2-3)FG+$$
$$4\tau^9(\tau^2-1)^2(102\tau^6-68\tau^4+21\tau^2-3)G^2=$$
$$-12\tau^9(\tau^2-1)(68\tau^6-45\tau^4+12\tau^2-1)F^2-$$
$$8\tau^9(\tau^2-1)(157\tau^6-106\tau^4+30\tau^2-3)FG+$$
$$4\tau^9(\tau^2-1)(102\tau^8-170\tau^6+89\tau^4-24\tau^2+3)G^2,$$
\begin{equation}\label{eq:7wkpfg33c}
(\vec w_{F,G,3}(\nu))_3=3\tau^4(\tau+1)(2\tau+1)^3F^2+
\end{equation}
$$t^4(t+1)(2t+1)(10t^2+10t+3)FG
-2\tau^5(\tau+1)^2(2\tau+1)((\tau+1)^3-\tau^3)G^2=$$
$$\tau^4(\tau+1)(2\tau+1)\times$$
$$3(2\tau+1)^2F^2+(10\tau^2+10\tau+3)FG-2(\tau+1)(3\tau^3+3\tau^2+\tau)G^2=$$
$$\tau^4(\tau+1)(2\tau+1)\times$$
$$3(4\tau^2+4\tau+1)F^2+(10\tau^2+10\tau+3)FG
-2(3\tau^4+6\tau^3+4\tau^2+\tau)G^2.$$
%%%%%%%%%%%%%%%%%%%%%%%%%%%%%%%%%%%%%%%%%%%%%%%%%%%%%%%%%%%%%%%%%%%%%%%%%%%%%
%%%%%%%%F\delta f_{1,0,k}(1,\nu+1)+G\delta^2f_{1,0,k}(1,\nu+1)%%%%%%%%%%%
%%%%%%%%%%%%%%%%%%%%%%%%%%%%%%%%%%%%%%%%%%%%%%%%%%%%%%%%%%%%%%%%%%%%%%%%%%%%
According to  (\ref{eq:yfg}), (\ref{eq:7.10}), (\ref{eq:7.14}),
 (\ref{eq:7wkpfg31c}), (\ref{eq:7wkpfg32c}), (\ref{eq:7wkpfg33c}),
\begin{equation}\label{eq:7.17h}
-\tau^4(\tau+1)^5w_{F,G,3,1}(\nu)y^{\ast\ast}_{F,G,k}(\nu+1)+
\end{equation}
$$w_{F,G,3,2}(\nu)y^{\ast\ast}_{F,G,k}(\nu)+
\tau^4(\tau-1)^5w_{F,G,3,3}(\nu)y^{\ast\ast}_{F,G,k}(\nu-1)=0.$$
%%%%%%%%%%%%%%%%%%%%%%%%%%%%%%%%%%%%%%%%%%%%%%%%%%%%%%%%%%%%%%%%%%%%%%%%%%
Since
$$f_{1,0,k}(1,\nu)=f_{1,0,k}^\ast(1,\nu)/(\nu+1)^2,$$
it follows from (\ref{eq:7.17h}), (\ref{eq:cfg2n}) -- (\ref{eq:cfg0n}) that
%%%%%%%%%%%%%%%%%%%%%%%%%%%%%%%%%%%%%%%%%%%%%%%%%%%%%%%%%%%%%%%%%%%%%%%%%%%
%%%%%%%%%%%%%%%%%%%%%%%%%%%%%%%%%%%%%%%%%%%%%%%%%%%%%%%%%%%%%%%%%%%%%%%%%%
\begin{equation}\label{eq:7.17ha}
c_{F,G,2}(\nu)x(\nu+1)+c_{F,G,1}(\nu)x(\nu)+
c_{F,G,0}(\nu)x(\nu-1)=0
\end{equation}
 for $x(\nu)=x_{F,G,k}(\nu),$
where
\begin{equation}\label{eq:xfgkn}
x_{F,G,k}(\nu)=F\delta f^\ast_k(1,\nu)+G\delta^2f^\ast_k(1,\nu)),\,
k=1,\,3.
\end{equation} Let
\begin{equation}\label{eq:bfgizn}
\beta_{F,G,i}^{\ast\ast}(z;\nu):=
F\beta^{\ast(1)}_{2i}(z;\nu)+G\beta_{2i}^{\ast(2)}(z;\nu)
\end{equation}
for $i=1,2.$ In view of (\ref{eq:15.4.hr8}),
$$\beta_{F,G,1}^{\ast\ast}(1;\nu):=
F\beta^{\ast(1)}_2(1;\nu)+G\beta_2^{\ast(2)}(1;\nu)=$$
$$F\delta f_1^\ast(1,\nu)+G\delta^2 f_1^\ast(1,\nu)=
x_{F,G,1}(\nu)$$
In view of (\ref{eq:15.4.hr1f}) with $j=1$ and (\ref{eq:xfgkn}),
\begin{equation}\label{eq:xfg3n}
x_{F,G,3}(\nu)=2\zeta(3)\beta_{F,G,1}^{\ast\ast}(1;\nu)
-\beta_{F,G,2}^{\ast\ast}(1;\nu)
\end{equation}
%%%%%%%%%%%%%%%%%%%%%%%%%%%%%%%%%%%%%%%%%%%%%%%%%%%%%%%%%%%
%%%%%%%%%%%%%%%%%%%%%%%%%%%%%%%%%%%%%%%%%%%%%%%%%%%%%%%%%%%%%%%%%%%%%%%%%%
Before to complete the proof of Theorem B, we want to check equality
 (\ref{eq:7.17ha}) for $\nu=1,\,k=3.$
%%%%%%%%%%%%%%%%%%%%%%%%%%%%%%check of 7.15a%%%%%%%%%%%%%%%%%%%%%%%%%%%%%
%%%%%%%%%%%%%%%%%%%%%%%%%%%%%%%%%%%%%%%%%%%%%%%%%%%%%%%%%%%%%%%%%%%%%%%%%%%
In view of (\ref{eq:xfg3n}), to check the equation (\ref{eq:7.17ha})
for $\nu=1$ it is  sufficient to check the equalities
$$%\begin{equation}\label{eq:7.17han1}
c_{F,G,2}(1)\beta_{F,G,i}^{\ast\ast}(2)+
c_{F,G,1}(1)\beta_{F,G,i}^{\ast\ast}(1)+
c_{F,G,0}(1)\beta_{F,G,i}^{\ast\ast}(0)=0,
$$%\end{equation}

In view of (\ref{eq:7wkpfg31c}) -- (\ref{eq:7wkpfg33c}),
$$%\begin{equation}\label{eq:cfg2n1}
c_{F,G,2}(1)=-54(-27F^2-23FG+28G^2),
$$%\end{equation}
$$%\begin{equation}\label{eq:eq:cfg1n1}
c_{F,G,1}(1)=12(-3679F^2-5646FG+5459G^2),
$$%\end{equation}
$$%\begin{equation}\label{eq:eq:cfg0n1}
c_{F,G,0}(1)=10(75F^2+63FG-228G^2)=30(25F^2+21FG-76F^2),
$$%\end{equation}
in view of (\ref{eq:a1r1i2n0}) -- (\ref{eq:a1r2i4n2}),
$$\beta^{\ast(1)}_{1,0,2}(1;2)=1029,\,\beta^{\ast(1)}_{1,0,2}(1;1)=34,\,
\beta^{\ast(1)}_{1,0,2}(1;0)=1,$$
$$\beta^{\ast(1)}_{1,0,4}(1;2)=\frac{14843}6,\,
\beta^{\ast(1)}_{1,0,4}(1;1)=\frac{327}4,\,
\beta^{\ast(1)}_{1,0,4}(1;0)=3.$$

$$\beta^{\ast(2)}_{1,0,2}(1;2)=2277,\,\beta^{\ast(2)}_{1,0,2}(1;1)=52,\,
\beta^{\ast(2)}_{1,0,2}(1;0)=1,$$
$$\beta^{\ast(2)}_{1,0,4}(1;2)=\frac{32845}6,\,
\beta^{\ast(2)}_{1,0,4}(1;1)=125,\,\beta^{\ast(2)}_{1,0,4}(1;0)=2,$$
$$\beta_{F,G,1}^{\ast\ast}(1,2)=1029F+2277G,$$
\begin{equation}\label{eq:betafg1}
\beta_{F,G,1}^{\ast\ast}(1,1)=
34F+52G,\,
\beta_{F,G,1}^{\ast\ast}(1,0)=F+G,
\end{equation}
$$\beta_{F,G,2}^{\ast\ast}(1,2)=(14843F+32845G)/6,$$
\begin{equation}\label{eq:betafg2}
\beta_{F,G,2}^{\ast\ast}(1,1)=(327F+500G)/4,
\beta_{F,G,2}^{\ast\ast}(1,0)=3F+2G.\end{equation}
Let further
$$\Delta_{i,\nu}(F,G)=c_{F,G,2}(\nu)\beta_{F,G,i}{\ast\ast}(\nu+1)+$$
$$c_{F,G,1}(\nu)\beta_{F,G,i}{\ast\ast}(\nu)+
c_{F,G,0}(\nu)\beta_{F,G,i}{\ast\ast}(\nu-1)$$
for $i=1,\,2\,\nu\in{\mathbb N}.$
We check now that $\Delta_{i,1}(F,G)=0$ for $i=1,2.$
We note that $\Delta_{i,1}(F,G)$ is homogenous polynomial relatively
variables $F,G$ of degree equal to $3.$ To establish the equalities
$%\begin{equation}\label{eq:di1FG0}
\Delta_{i,1}(F,G)=0
$%\end{equation}
with $i=1,\,2$
it is sufficient to check them in four points
$$(F,G)=(0,1),\,(1,0),\,(1,1),\,(1,2).$$
We have
$$c_{0,1,2}(1)=-54\times28=-12\times126,\,c_{0,1,1}(1)=12\times5521,$$
$$c_{0,1,0}(1)=-12\times190,$$
$$\beta_{0,1,1}^{\ast\ast}(2)=2277,\,\beta_{0,1,1}^{\ast\ast}(1)=52,\,
\beta_{0,1,1}^{\ast\ast}(0)=1,$$
$$\beta_{0,1,2}^{\ast\ast}(2)=32845/6,\,\beta_{0,1,2}^{\ast\ast}(1)=125,\,
\beta_{0,1,2}^{\ast\ast}(0)=2,$$
$$\Delta_{1,1}(0,1)=12\times(-126\times2277+5521\times52-190)=$$
$$12(-287092+286902-190)=0,$$
$$\Delta_{2,1}(0,1)=12\times(-21\times32845+5521\times125-190\time2)=$$
$$60(-21\times6569+5521\times25-76)=60(-137949+138025-76)=0,$$
$$c_{1,0,2}(1)=-54\times(-27)=6\times243,\,c_{1,0,1}(1)=
-12\times3679=-4\times11037,$$
$$c_{1,0,0}(1)=10\times75,$$
$$\beta_{1,0,1}^{\ast\ast}(2)=1029,\,\beta_{1,0,1}^{\ast\ast}(1)=34,\,
\beta_{1,0,1}^{\ast\ast}(0)=1,$$
$$\beta_{1,0,2}(2)=14843/6,\,\beta_{1,0,2}(1)=327/4,\,\beta_{1,0,2}(0)=3,$$
$$\Delta_{1,1}(1,0)=6(243\times1029-7358\times34+125)=$$
$$6(250047-250172+125)=0.$$
$$\Delta_{2,1}(1,0)=243\times14843-11037\times327+2250=$$
$$3606849-3609099+2250=3609099-3609099=0.$$
$$c_{1,1,2}(1)=-54\times(-22)=36\times33,\,c_{1,1,1}(1)=
-12\times3804=-36\times1268,$$
$$c_{1,1,0}(1)=-900=-36\times25,$$
$$\beta_{1,1,1}^{\ast\ast}(2)=2\times1653,\,\beta_{1,1,1}^{\ast\ast}(1)=
2\times43,\,
\beta_{1,1,1}^{\ast\ast}(0)=2,$$
$$\beta_{1,1,2}^{\ast\ast}(2)=7948,\,\beta_{1,1,2}^{\ast\ast}(1)=827/4,\,
\beta_{1,1,2}^{\ast\ast}(0)=5,$$
$$\Delta_{1,1}(1,1)=72(33\times1653-1268\times43-25)=$$
$$72(54559-54524-25)=0,$$
$$\Delta_{2,1}(1,1)=36(33\times7948-317\times827-125)=$$
$$36(262284-262159-125)=0.$$
$$c_{1,2,2}(1)=-54\times(39)=-18\times117\,c_{1,2,1}(1)=
18\times4742,$$
$$c_{1,1,0}(1)=-18\times395,$$
$$\beta_{1,2,1}^{\ast\ast}(2)=3\times1861,\,
\beta_{1,2,1}^{\ast\ast}(1)=3\times46,\,
\beta_{1,2,1}^{\ast\ast}(0)=3,$$
$$\beta_{1,2,2}^{\ast\ast}(2)=80533/6,\,\beta_{1,2,2}^{\ast\ast}(1)=1327/4,\,
\beta_{1,2,1}^{\ast\ast}(0)=7,$$
$$\Delta_{1,1}(1,2)=54(-117\times1861+4742\times46-395)=$$
$$54(-217737+218132-395)=0.$$
$$\Delta_{2,1}(1,2)=-351\times80533+9times2371\times1327-63\times790)=$$
$$9(-39\times80533+271\times1327-7\times790)=$$
$$(-3140787+3146317-5530=0.$$
So  $\Delta_{i,1}(F,G)=0$ for any $F$ and $G.$ Therefore
the equality (\ref{eq:7.17ha}) holds

\noindent for $\nu=1$ and any $\{F,G\}\subset{\mathbb R}.$
%%%%%%%%%%%%%%%%%%%%%%%%%%%%%%%%%%%%%%%%%%%%%%%%%%%%%%%%%%%%%%%%%%%%
\refstepcounter{section}
{\begin{center}\large\bf\S 6. Auxilliary continued fraction.
\end{center}}
%%%%%%%%%%%%%%%%%%%%%%%%%%%%%%%%%%%%%%%%%%%%%%%%%%%%%%%%%%%%%%%%%%%%%%%
Let
\setcounter{equation}0
\begin{equation}\label{eq:cfgknast}
c_{u,v,k}^\ast(\nu)=c_{u,v,k}(\nu)/(u+v)^2,
\end{equation}
for $k=0,\,1,\,2,$
\begin{equation}\label{eq:bfgnp1ast}
b_{u,v}^\ast(\nu+1)=-c_{u,v,1}^\ast(\nu)\in{\mathbb Q}[u,v]/(u+v)^2
\end{equation}
for $\nu\in\mathbb N,$
\begin{equation}\label{eq:afgnp1ast}
a_{u,v}^\ast(\nu+1)=-c_{u,v,0}^\ast(\nu)c_{u,v,2}^\ast(\nu-1)
\end{equation}
for $\nu\in[2+\infty)\cap{\mathbb N},$
\begin{equation}\label{eq:avg2ast}
a_{u,v}^\ast(2)=-c_{u,v,0}^\ast(1),
\end{equation}
\begin{equation}\label{eq:bfg0ast}
P^\ast_{u,v}(0)=b_{u,v}^\ast(0)=(3u+2v)/(u+v),\,Q_{u,v}(0)=1,
\end{equation}
\begin{equation}\label{eq:bfg1ast}
Q_{u,v}^\ast(1)=b_{u,v}^\ast(1)=(34u+52v)(u+v)
\end{equation}
\begin{equation}\label{eq:Pfg1ast}
P_{u,v}(1)=(327u+500v)(4u+4v)
\end{equation}
\begin{equation}\label{eq:afg1ast}
a_{u,v}^\ast(1)=P_{u,v}^\ast(1)=-b_{u,v}^\ast(0)b_{u,v}^\ast(1).
\end{equation}
  Let us consider the continued fraction,
\begin{equation}\label{eq:G0ast}
b_{u,v}(0)^\ast+\frac{a_{u,v}^\ast(1)\vert}{b_{u,v}^\ast(1)}+
\frac{a_{u,v}^\ast(2)\vert}{b_{u,v}^\ast(2)}
+\frac{a_{u,v}^\ast(3)\vert}{b_{u,v}^\ast(3)}+
\frac{a_{u,v}^\ast(4)\vert}{\vert b_{F,G}^\ast(4)}... .
\end{equation}
Let $r_{u,v}^\ast(\nu)$ be the $\nu$-th convergent of this continued fraction.
Let $P_{u,v}^\ast(\nu)$ and $Q_{u,v}^\ast(\nu)$ be respectively nominator
and denominator of convergent $r_{u,v}^\ast(\nu).$
%%%%%%%%%%%%%%%%%%%%%%%%%%%%%%%%%%%%%%%%%%%%%%%%%%%%%%%%%%%%%%%%%%%%%%%%
Let us consider the equations
\begin{equation}\label{eq:7.17hhh}
c_{F,G,2}(\nu)x_{\nu+1}+c_{F,G,1}(\nu)x_{\nu}+c_{F,G,0}(\nu)x_{\nu-1}=0,
\end{equation}
where $\nu\in{\mathbb N}.$ If $F+G\ne0,$ then the equation $\ref{eq:7.17hhh}$
is equivalent to the equation
\begin{equation}\label{eq:7.17hhhast}
c_{F,G,2}^ast(\nu)x_{\nu+1}+c_{F,G,1}^\ast(\nu)x_{\nu}+
c_{F,G,0}^\ast(\nu)x_{\nu-1}=0,
\end{equation}
 It follows from (\ref{eq:7.17ha}) that
 $$x_\nu=x_{F,G,k}(\nu)=
F\delta f^\ast_{1,0,k}(1,\nu)+G\delta^2f_{1,0,k}(1,\nu)$$
satisfies to the equation (\ref{eq:7.17hhh}) for $\nu\in{\mathbb N}$ and
 fixed $k\in\{1,\,3\}.$

If $G\ne0,$ then, in view of (\ref{eq:cfg2n}) -- (\ref{eq:cfg0n}),
$$c_{F,G,2}(\nu)=-12\tau^8G^2(1+o(1)) (\tau\to\infty),$$
$$c_{F,G,1}(\nu)=408\tau^8G^2(1+o(1)) (\tau\to\infty),$$
$$c_{F,G,0}(\nu)=-12\tau^8G^2(1+o(1)) (\tau\to\infty).$$
If $G=0,$ then, in view of (\ref{eq:cfg2n}) -- (\ref{eq:cfg0n}),
$$c_{F,G,2}(\nu)=24\tau^6F^2(1+o(1)) (\tau\to\infty),$$
$$c_{F,G,1}(\nu)=-24\times64\tau^6F^2(1+o(1)) (\tau\to\infty),$$
$$c_{F,G,0}(\nu)=24\tau^6F^2(1+o(1)) (\tau\to\infty).$$
In any case the equation (\ref{eq:7.17hhh}) is difference equation
of Poincar\'e type  with characteristic polynomial  $\lambda^2-34\lambda+1.$
 Hence, if $\{x_\nu\}_{\nu=1}^{+\infty}$ is a non-zero solution
  of (\ref{eq:7.17hhh}), $\varepsilon\in(0,1),$
   then there are $C_1(\varepsilon)>0$ and $C_2(\varepsilon)>0$ such that
 only two possibilities exist:
 \begin{equation}\label{eq:7.bf1}
C_1(\varepsilon)\left(1+\sqrt{2}\right)^{-4\nu(1+\varepsilon)}\le
\vert x_\nu\vert\le
C_2(\varepsilon)\left(1+\sqrt{2}\right)^{-4\nu(1-\varepsilon)}
\end{equation}
for all $\nu\in{\mathbb N}$ or
 \begin{equation}\label{eq:7.bf2}
C_1(\varepsilon)\left(1+\sqrt{2}\right)^{4\nu(1-\varepsilon)}\le
\vert x_\nu\vert\le
C_2(\varepsilon)\left(1+\sqrt{2}\right)^{4\nu(1+\varepsilon)}.
\end{equation}
for all $\nu\in{\mathbb N}.$
In view of (\ref{eq:15.4.hr8}), if
 $x_\nu=\beta^{\ast(r)}(1;\nu)=\delta^r f_1^\ast(1,\nu)$
with $r=0,1,2,$
% or $x_\nu=\delta f_{1,0,1}^\ast(1,\nu)+2\delta^2 f_{1,0,1}^\ast(1,\nu),$
%$x_\nu=F\delta f_{1,0,1}^\ast(1,\nu)+G\delta^2 f_{1,0,1}^\ast(1,\nu),$
 then (\ref{eq:7.bf1}) is impossible. Therefore
 \begin{equation}\label{eq:betaastn_1}
C_1(\varepsilon)\left(1+\sqrt{2}\right)^{4\nu(1-\varepsilon)}\le
 \beta^{\ast(r)}(1;\nu) \le
C_2(\varepsilon)\left(1+\sqrt{2}\right)^{4\nu(1+\varepsilon)}.
\end{equation}
for $r=0,1,2,$ and all $\nu\in{\mathbb N}.$

\bf Lemma 6.1. \it The following equalities hold:

\begin{equation}\label{eq:l6.1}
\lim\limits_{\nu\to\infty}\beta^{\ast(1)}_{1,0,2}(1;\nu))=+\infty,
\lim\limits_{\nu\to\infty}
\beta^{\ast(2)}_2(1;\nu)/\beta^{\ast(1)}_2(1;\nu)=+\infty,
\end{equation}

\bf Proof. \rm The first of equalities (\ref{eq:l6.1}) is obvious.
We prove the second of equalities (\ref{eq:l6.1}). \rm
According to Stirling's formula:
$$\log(n!)=(n+1/2)\log(n+1)-n+O(1).$$
Let $\beta> 0,\,n=\beta\nu+\eta_{1,\nu}\in{\mathbb Z},$
\begin{equation}\label{eq:logn!c}
\nu\in[(2(\beta+1)/\beta,+\infty)\cap{\mathbb Z},\,\vert \eta_\nu\vert<2.
\end{equation}
Then (for fixed $\beta$)
\begin{equation}\label{eq:logn!r}
\log(n!)=(\beta\nu+\eta_\nu+1/2)\log(\beta\nu+\eta_\nu+1)-\beta\nu+O(1)=
\end{equation}
$$(\beta\nu+\eta_\nu+1/2)\log(\beta\nu)-\beta\nu+O(1)=
(\beta\nu)\log(\beta\nu)-\beta\nu+O(1)\log(\nu).$$
Clearly,
$$\log((\nu+1)!)=\nu\log(\nu)-\nu+O(1)\log(\nu+1).$$

Let further $\gamma\in(0,1),\,k=[\gamma\nu],$
where $\nu\in[2(\gamma+1)/\gamma,+\infty)\cap{\mathbb Z}.$
Then, in view of (\ref{eq:logn!c}) -- (\ref{eq:logn!r})
with $\beta=\gamma$ and $k$ in the role of $n,$
\begin{equation}\label{eq:logk!}
\log(k!)=\gamma\nu\log(\gamma\nu)-\gamma\nu+O(1)\log(\nu).
\end{equation}
Further we have
 $\nu+1-k=\nu+1-\gamma\nu+\{\gamma\nu\}=(1-\gamma)\nu+1+\{\gamma\nu\}.$

\noindent If $\nu\in[2(2-\gamma)/(1-\gamma),+\infty)\cap{\mathbb Z},$ then,
 according to  (\ref{eq:logn!c}) -- (\ref{eq:logn!r})

\noindent with $\beta=1-\gamma,$ and $n+1-k$ in the role of $n,$
we have the equality
\begin{equation}\label{eq:lognp1mk!}
\log((\nu+1-k)!)=(1-\gamma)\nu\log((1-\gamma)\nu)-(1-\gamma)\nu+O(1)\log(\nu).
\end{equation}

Since $\nu+k=\nu+\gamma\nu-\{\gamma\nu\}=(1+\gamma)\nu-\{\gamma\nu\}$ it
follows that,

\noindent if $\nu\in[2(2+\gamma)/(1+\gamma),+\infty)\cap{\mathbb Z},$ then,
 in view of (\ref{eq:logn!c}) -- (\ref{eq:logn!r})

\noindent with $\beta=1+\gamma$ and $n+k$ in the role of $n,$
we have the equality
\begin{equation}\label{eq:lognpk!}
\log((\nu+k)!)=(1+\gamma)\nu\log((1+\gamma)\nu)-(1+\gamma)\nu+O(1)\log(\nu).
\end{equation}
Let $\nu\in[2(1/\min(\gamma,1-\gamma)+1).$ Then
(\ref{eq:logk!}) -- (\ref{eq:lognpk!}) hold, and,
moreover,
$$\log((\nu+1)!)=\nu\log(\nu)-\nu+O(1)\log(\nu),$$
$$\log(\nu!)=\nu\log(\nu)-\nu+O(1)\log(\nu).$$
So, if $\nu\in[2/\min(\gamma,1-\gamma)+2,+\infty),k=[\gamma\nu],$
then
$$
\log\left(\binom{\nu+1}k\right)=\nu\log(\nu)-\nu-(\gamma\nu\log(\gamma\nu)
-\gamma\nu)-
$$
$$((1-\gamma)\nu\log((1-\gamma)\nu)-(1-\gamma)\nu)+O(1)\log(\nu)=$$
$$\nu\log(\nu)-(\gamma\nu\log(\gamma\nu)-(\gamma\nu\log(\gamma)-$$
$$((1-\gamma)\nu\log(\nu)-((1-\gamma)\nu\log(1-\gamma)+O(1)\log(\nu)=$$
$$(\gamma\log(1/\gamma)+(1-\gamma)\log(1/(1-\gamma)))\nu+O(1)\log(\nu),$$
and, analogously,
$$\log\left(\binom{\nu+k}k\right)=
((1+\gamma)\log(1+\gamma)+\gamma\log(1/\gamma))\nu
+O(1)\log(\nu)$$
Therefore, if $\nu\in[2/\min(\gamma,1-\gamma)+2,+\infty),k=[\gamma\nu],$ then
\begin{equation}\label{eq:beta02gnn}
\log\left(\binom{\nu+1}k\binom{\nu+k}k\right)=\psi_1(\gamma)\nu+O(1)\log(\nu),
\end{equation}
where $\psi_1(\gamma)=(1+\gamma)\log(1+\gamma)-(1-\gamma)\log(1-\gamma)-
2\gamma)\log(\gamma).$ Clearly,
$$\left(\frac d{d \gamma} \psi_1\right)(\gamma)=
\log((1-\gamma^2)/\gamma^2)),$$
and
$$\psi_1(\gamma)<\psi_1\left(1/\sqrt{2}\right)=$$
$$\log\left(\left(1+1/\sqrt{2}\right)/\left(1-1/\sqrt{2}\right)\right)+$$
$$\left(1/\sqrt{2}\right)
\log\left(\left(1+1/\sqrt{2}\right)\left(1-1/\sqrt{2}\right)\right)-$$
$$
\left(1/\sqrt{2}\log(1/2)\right)=2\log\left(\left(1+\sqrt{2}\right)\right),
$$
if $\gamma\ne1/\sqrt{2}.$
We can rewrite (\ref{eq:betaastn_1}) in the form
 \begin{equation}\label{eq:betaastn_2}
C_1(\varepsilon)\exp
\left(2\psi_1\left(1/\sqrt{2}\right)\nu(1-\varepsilon)\right)\le
 \beta^{\ast(r)}(1;\nu) \le\end{equation}
$$
C_2(\varepsilon)\exp
\left(2\psi_1\left(1/\sqrt{2}\right)\nu(1+\varepsilon)\right)
$$
for $r=0,1,2,$ and all $\nu\in{\mathbb N}.$

In view of (\ref{eq:6bi2}), let
$$\psi_{2,\nu}(k)=\beta_{2,k+1,\nu}^{(0)}/\beta_{2,k,\nu}^{(0)}=
 (\nu+k+1)(\nu-k+1)/(k+1)^2.$$
Then $\psi_{2,\nu}(k)=1$ if and only if
$2k^2+2k-\nu(\nu+2)=0.$ Let
$$r(\nu)=-1/2+\sqrt{1/4+\nu(\nu+2)/2}=(1+O(1)/\nu)\nu/\sqrt{2}$$
Clearly, $\beta_{2,k,\nu}^{(0)}$ increases together with increasing of
 $k\in[0,r(\nu)]\cap{\mathbb Z}$ and
$\beta_{2,k,\nu}^{(0)}$ decreases together with increasing of
 $k\in(r(\nu),\nu+1]\cap{\mathbb Z}.$

We fix $\gamma_1\in\left(0,1/\sqrt{2}\right)$ and
$\gamma_2\in\left(1/\sqrt{2},1\right).$
 Clearly, there exists $\nu_0\in{\mathbb N}$ such that
$$\gamma_1<r(\nu)/\nu<\gamma_2$$ for
all $\nu\in[\nu_0,+\infty)\cap{\mathbb N}.$
Let $\nu>\max(\nu_0,1/\gamma_1,1/(1-\gamma_2)).$ Then, in view

\noindent of (\ref {eq:beta02gnn}) we have
$$\sum\limits_{0\le k\le\gamma_1\nu}\beta^{(r)}_{2,k,\nu}=
\sum\limits_{0\le k\le\gamma_1\nu}k^r\beta^{(0)}_{2,k,\nu}<
[\gamma_1\nu]^{r+1}\beta^{(0)}_{2,[\gamma_1\nu],\nu}=$$
$$\exp(2\psi_1(\gamma_1)\nu+O(1)\log(\nu)),$$
and analogously
$$\sum\limits_{\gamma_2\nu\le k\le\nu+1}\beta^{(r)}_{2,k,\nu}=
\exp(2\psi_1(\gamma_2)\nu+O(1)\log(\nu)).$$
Let
$$\gamma_3=\min(\psi_1(1/\sqrt{2})-\psi_1(\gamma_1),
\psi_1(1/\sqrt{2})-\psi_1(\gamma_2)).$$
Then, in view of (\ref{eq:betaastn_2})
$$\beta^{(r)}_2(1,\nu)=
\left(\sum\limits_{\gamma_1\nu<k<\gamma_2\nu}k^r\beta^{(0)}_{2,k,\nu}\right)
\times$$
$$(1+O(1)\exp(-2\gamma_3\nu+O(\log(\nu)).$$
Hence, there exists $\nu_1\in{\mathbb N}$ such that
for all $\nu\in[\nu_1,+\infty)cap{\mathbb Z}$ we have
$$\beta^{(2)}_2(1,\nu)=
\left(\sum\limits_{\gamma_1\nu<k<\gamma_2\nu}k^2\beta^{(0)}_{2,k,\nu}\right)
\times$$
$$(1+O(1)\exp(-2\gamma_3\nu+O(\log(\nu))\ge$$
$$[\gamma_1\nu]
\left(\sum\limits_{\gamma_1\nu<k<\gamma_2\nu}k\beta^{(0)}_{2,k,\nu}\right)
\times$$
$$(1+O(1)\exp(-2\gamma_3\nu+O(\log(\nu)),$$ and
$$\beta^{(1)}_2(1,\nu)=
\left(\sum\limits_{\gamma_1\nu<k<\gamma_2\nu}k\beta^{(0)}_{2,k,\nu}\right)
(1+O(1)\exp(-2\gamma_3\nu+O(\log(\nu)).$$ So,
$\beta^{(2)}_2(1,\nu)\ge[\gamma_1\nu]\beta^{(2)}_2(1,\nu)(1+o(1)),$ when
$\nu\to+\infty.$ $\square.$

 Let conditions  Theorem B are fulfilled. Then, in
view of (\ref{eq:barzn}),
\begin{equation}\label{eq:xfg1ne0}
\beta_{F,G,1}^{\ast\ast}(1;\nu)= (F+G)\beta^{\ast(1)}_{1,0,2}(1;\nu)+
G(\beta^{\ast(2)}_{1,0,2}(1;\nu)-\beta^{\ast(1)}_{1,0,2}(1;\nu))\ne0\end{equation}
for $\nu\in{\mathbb N}_0,$ and therefore, in view of (\ref{eq:l6.1}),
\begin{equation}\label{eq:xfg1nti}
\beta_{F,G,1}^{\ast}(1;\nu)= \beta^{\ast(1)}_2(1;\nu)
(F+G\beta^{\ast(2)}_2(1;\nu)/\beta^{\ast(1)}_2(1;\nu))\to\infty,
\end{equation}
when $\nu\to\infty.$ Moreover, if $F\ne0$ and $G/F\not\in{\mathfrak B}$ then
(\ref{eq:xfg1ne0}) and (\ref{eq:xfg1nti}) hold, and, if in this case
 $x_\nu=x_{F,G,1}(\nu)=\beta_{F,G,1}^{\ast\ast}(1;\nu),$
 then (\ref{eq:7.bf1}) is impossible.

In view of (\ref{eq:d}) with $\alpha=1,$ (\ref{eq:15.4.hr}) and
(\ref{eq:15.4.hr1f}) with $j=1,$
$$\delta^r f_{1,0,3}^\ast(1,\nu)=(\nu+1)^2O(1);$$ hence, if
$x_\nu=x_{F,G,3}(\nu)=F\delta f_{1,0,3}^\ast(1,\nu)+
G\delta^2 f_{1,0,3}^\ast(1,\nu),$
 then (\ref{eq:7.bf2}) is impossible. Therefore
 \begin{equation}\label{eq:7.bf3}
\frac{C_1(\varepsilon)/C_2(\varepsilon)}
{\left(1+\sqrt{2}\right)^{8\nu(1+\varepsilon)}}\le
\bigg\vert2\zeta(3)-\frac{\beta_{F,G,2}^{\ast\ast}(1,\nu)}
{\beta_{F,G,1}^{\ast\ast}(1,\nu)}\bigg\vert\le
\frac{C_2(\varepsilon)/C_1(\varepsilon)}
{\left(1+\sqrt{2}\right)^{-8\nu(1-\varepsilon)}}.
\end{equation}
%%%%%%%%%%%%%%%%%%%%%%%%%%%%%%%%%%%%%%%%%%%%%%%%%%%%%%%%%%%%%%%%%%%%%%%%
Let
\begin{equation}\label{eq:dfgnu}
\delta_{u,v}^\ast(\nu)=\prod\limits_{j=1}^{\nu-1}c_{u,v,2}^\ast(j)
\,\,\text{for}\,\,\nu\in{\mathbb N}_0.
\end{equation}
So, $\delta_{F,G}^\ast(1)=\delta_{F,G}^\ast(0)=1.$
Let $c_{F,G,1}^{\ast\ast}(\nu)=c_{F,G,1}^{\ast}(\nu)$
for all $\nu\in{\mathbb N},$ let
$$c_{F,G,0}^{\ast\ast}(\nu)=c_{F,G,0}^{\ast}(\nu)c_{F,G,0}^{\ast}(\nu-1)$$
for all $\nu\in[2,+\infty)\cap{\mathbb N},$ and let
$c_{F,G,0}^{\ast\ast}(1)=c_{F,G,0}^{\ast}(1).$
 If  conditions  Theorem B are fulfilled, then
$\delta_{F,G}(\nu)\ne0$ far all the $\nu\in{\mathbb N},$ and
  the equation (\ref{eq:7.17hhh}) and the system
$$%\begin{equation}\label{eq:7.18hhh}
\bigg\{\begin{matrix}
 y_{\nu+1}+c_{F,G,1}^{\ast}(\nu)y(\nu)+
c_{F,G,0}^\ast(\nu)c_{F,G,2}^ast(\nu-1)y(\nu-1)=0\\
 y_\nu=\delta_{F,G}(\nu)^\ast x_\nu,\,\nu\in{\mathbb N}
\end{matrix}$$
are equivalent. Moreover,
$P_{F,G}^\ast(\nu)$ and
$\delta_{F,G}^\ast(\nu)\beta_{F,G,2}^{\ast\ast}(1,\nu)/
\beta_{F,G,1}^{\ast\ast}(1,0)$
satisfy to the first of equations (\ref{eq:7.17hhh})
and the same initial conditions. Therefore
\begin{equation}\label{eq:Pbeta}
P_{F,G}^\ast(\nu)=\delta_{F,G}(\nu)\beta_{F,G,2}^{\ast\ast}(1,\nu)/
\beta_{F,G,1}^{\ast\ast}(1,0),
\end{equation}
Analogously, we have
\begin{equation}\label{eq:Qbeta}
Q_{F,G}^\ast(\nu)=\delta_{F,G}(\nu)\beta_{F,G,1}^{\ast\ast}(1,\nu)/
\beta_{F,G,1}^{\ast\ast}(1,0),
\end{equation}
$
r^\ast_{F,G}(\nu)=\beta_{F,G,2}^{\ast\ast}(1,\nu)/
\beta_{F,G,1}^{\ast\ast}(1,\nu),
$
for all $\nu\in{\mathbb N}_0.$
In view of (\ref{eq:7.bf3})
$$\lim\limits_{\nu\to\infty} r_{F,G}^\ast(\nu)=2\zeta(3).$$
\refstepcounter{section}
{\begin{center}\large\bf\S 7. End  of the proof of Theorem B.
\end{center}}

\bf Lemma 7.1. \it The following equalities hold:
\setcounter{equation}0
\begin{equation}\label{eq:7.1h}
P^\ast_{u,v}(\nu)=P^\ast_{u/v,1}(\nu),\,
Q^\ast_{u,v}(\nu)=Q_{u/v,1}^\ast(\nu),\delta_{u,v}(\nu)=\delta d_{u/v,1}(\nu).
\end{equation}

\bf Proof. \rm In view of  (\ref{eq:Pfg1})
If $\nu=0,\,1,$ then the equalities (\ref{eq:7.1h}) directly
follows from (\ref{eq:bfg0ast}) -- (\ref{eq:Pfg1ast}) and (\ref{eq:dfgnu})
In view of (\ref{eq:cfgknast}),
\begin{equation}\label{eq:cfgknh}
c_{u,v,k}^\ast(\nu)=c_{u/v,1,k}^\ast(\nu)
\end{equation}
for $k=0,\,1,\,2,\nu\in{\mathbb N}.$ Therefore the last equality
in (\ref{eq:7.1h}) holds for all $\nu\in{\mathbb N}.$
In view (\ref{eq:afg1ast})), (\ref{eq:bfg0ast}), (\ref{eq:bfg1ast}),
 (\ref{eq:afgnp1ast}), (\ref{eq:bfgnp1ast}) and (\ref{eq:cfgknh}),
$$%\begin{equation}\label{eq:abuvn}
a_{u,v}(\nu)=a_{u/v,1}(\nu),\,b_{u,v}(\nu)=b_{u/v,1}(\nu)
$$%\end{equation}
for $\nu\in[2,+\infty)\cap{\mathbb Z}.$
Let $\nu\in[2,+\infty)\cap{\mathbb Z},$ and let (\ref{eq:7.1h}) hold
for all $\nu-\kappa$ with $\kappa\in[0,\nu-1]\cap{\mathbb Z}.$
Then we have
$$P_{u,v}^\ast(\nu)=b_{u,v}^\ast(\nu)P_{u,v}^\ast(\nu-1)+a_{u,v}^\ast(\nu)
P_{u,v}^(\nu-2)=$$
$$b_{u/v,1}^\ast(\nu)P_{u/v,1}^\ast(\nu-1)+
a_{u/v,1}^\ast(\nu)P_{u/v,1}^\ast(\nu-2)=P_{u/v,1}^\ast(\nu)$$
$$Q_{u,v}^\ast(\nu)=b_{u,v}^\ast(\nu)Q_{u,v}^\ast(\nu-1)+
a_{u,v}^\ast(\nu)Q_{u,v}^\ast(\nu-2)=$$
$$b_{u/v,1}^\ast(\nu)Q_{u/v,1}^\ast(\nu-1)+
a_{u/u,1}^\ast(\nu)Q_{u/v,1}^\ast(\nu-2)=Q_{u/v,1}(\nu).$$
$\square$

\bf Lemma 7.2. \it If conditions of the Theorem B
are fulfilled, then

\begin{equation}\label{eq:deltast}
\delta_{u,v}^\ast(\nu)=\delta_{u,v}(\nu)/(u+v)^{\max{2\nu-2,0}}
\end{equation}
where $\delta^\ast_{u,v}(\nu)$ is homogeneous polynomial
in ${\mathbb Z}[u,v],$
and
\begin{equation}\label{eq:degdeltast}
\max(2\nu-2,0)=\deg_u(\delta_{u,v}(\nu))=\deg_v(\delta_{u,v}(\nu))=
\deg(\delta_{u,v}(\nu)).
\end{equation}

\bf Proof. \rm We have to prove the last equality in (\ref{eq:degdeltast}),
because other assertions of the Lemma are obvious.
In view of (\ref{eq:cfg0n}) and (\ref{eq:cfg2n}),
$$c_{F,G,0}(\nu)=-frac{(\tau-1)^2(2\tau+1)}
{(\tau+1)^2(2\tau-1)}\times$$
$$\tau(\tau+1)c_{F,G,2}(\nu+1)\ne0$$
for all $\nu\in{\mathbb N}_0.$ Therefore the last equality
in (\ref{eq:degdeltast}) holds also.

In view  of (\ref{eq:betafg1}) -- (\ref{eq:betafg2}),
(\ref{eq:Pbeta}) -- (\ref{eq:Qbeta}) and (\ref{eq:deltast}),
the  following equalities hold:
$$
P_{u,v}(\nu)=4P_{u,v}^{\ast}(\nu)(u+v)^{2\nu}=
$$
$$16(u+v)\delta_{u,v}(\nu)\beta^{\ast\ast}_{u,v,2}(1,\nu),$$
$$
Q_{u,v}(\nu)=Q_{u,v}^{\ast}(\nu)(u+v)^{2\nu-1}=
$$
$$4\delta_{u,v}(\nu)\beta{\ast\ast}_{u,v,1}(1,\nu)$$
$$\max(2\nu,1)=\deg_u(P_{u,v}(\nu))=\deg_v(P_{u,v}(\nu))=\deg(P_{u,v}(\nu)),$$
$$\max(2\nu-1,0)=\deg_u(Q_{u,v}(\nu))=\deg_v(Q_{u,v}(\nu))=\deg(Q_{u,v}(\nu)),$$
where $\nu\in {\mathbb N}_0.$ $\square$
%%%%%%%%%%%%%%%%%%%%%%%%%%%%%%%%%%%%%%%%%%%%%%%%%%%%%%%%%%%%%%%%%%%
%%%%%%%%%%%%%%%%%%%%%%% References %%%%%%%%%%%%%%%%%%%%%%%%%

%%%%%%%%%%%%%%%%%%%%%% ADDRESS %%%%%%%%%%  ADDRESS %%%%%%%%%%%%%%%%%%%%%%
\vskip 10pt
 {\it E-mails:}
$$l_-a_-gutnik33\text{$@$}mail.ru,$$
\hspace*{4cm}\ \ \ \ \ \ \ \ {\sl\ gutnik$@$gutnik.mccme.ru}.

%%%%%%%%%%%%%%%%%%%%%%%%%%%%%%%%%%%%%%%%%%%%%%%%%%%%%%%%%%%%%%%%%%%%%%%%%

\begin{thebibliography}{19}
\bibitem[R.Ap\'ery(1981)R.Ap\'ery]{Apery}
R.Ap\'ery, Interpolation des fractions continues et irrationalite de certaines
constantes, Bulletin de la section des sciences du C.T.H., 1981, No 3,
 37 -- 53.
 \bibitem[Oskar Perron(1996)Oskar Perron]{Perron}
Oskar Perron, Die Lehre von den Kettenbr\"uche. Dritte, verbesserte und
erweiterte Auflage. 1954 B.G.Teubner Verlaggesellshaft. Stuttgart.
 \bibitem[Yu.V. Nesterenko(1996)Yu.V. Nesterenko]{Nest}
Yu.V. Nesterenko, A Few Remarks on $\zeta(3),$ Mathematical Notes, Vol 59, No
6, 1996,
 Matematicheskie Zametki, 1996, Vol 59, No 6, pp. 865 -- 880,(in Russian).
 \bibitem[L.A.Gutnik,1(2002)L.A.Gutnik]{G1}
 L.A.Gutnik, On linear forms with coefficients
in ${\mathbb N}\zeta(1+{\mathbb N})$ (the detailed version,part 3),
 Max-Plank-Institut f\"ur Mathematik, Bonn,
 Preprint Series, 2002, 57, 1 -- 33.
 \bibitem[L.A.Gutnik,2(2005)L.A.Gutnik]{G2}
 L.A.Gutnik, On the measure of nondiscreteness
of some modules, Max-Plank-Institut f\"ur Mathematik, Bonn, Preprint Series,
2005, 32, 1 -- 51.
\bibitem[L.A.Gutnik,3(2006)L.A.Gutnik]{G3}
L.A.Gutnik, On the Diophantine approximations of logarithms in cylotomic
fields,
 Max-Plank-Institut f\"ur Mathematik, Bonn,
Preprint Series, 2006, 147, 1 -- 36.
 \bibitem[L.A.Gutnik,4(2006)L.A.Gutnik]{G4}
L.A.Gutnik, On some systems of difference equations, part 1,
Max-Plank-Institut f\"ur Mathematik, Bonn,
 Preprint Series, 2006, 23, 1 -- 37.
\bibitem[L.A.Gutnik,5(2006)L.A.Gutnik]{G5}
L.A.Gutnik,, On some systems of difference equations, part 2,
Max-Plank-Institut f\"ur Mathematik, Bonn, Preprint Series,2006, 49, 1 -- 31.
 \bibitem[L.A.Gutnik,6(2006)L.A.Gutnik]{G6}
L.A.Gutnik, On some systems of difference equations, part 3,
Max-Plank-Institut f\"ur Mathematik, Bonn,
 Preprint Series, 2006, 91, 1 -- 52.
 \bibitem[L.A.Gutnik,7(2006)L.A.Gutnik]{G7}
L.A.Gutnik, On some systems of difference equations, part 4,
Max-Plank-Institut f\"ur Mathematik, Bonn,
 Preprint Series, 2006, 101, 1 -- 49.
 \bibitem[L.A.Gutnik,8(2006)L.A.Gutnik]{G8}
L.A.Gutnik, On some systems of difference equations, part 5,
Max-Plank-Institut f\"ur Mathematik, Bonn,
 Preprint Series, 2006, 115, 1 -- 9.
 \bibitem[L.A.Gutnik,9(2006)L.A.Gutnik]{G9}
L.A.Gutnik, On some systems of difference equations, part 6,
Max-Plank-Institut f\"ur Mathematik, Bonn,
 Preprint Series, 2007, 16, 1 -- 30.
 \bibitem[L.A.Gutnik,10(2007)L.A.Gutnik]{G10}
L.A.Gutnik, On some systems of difference equations, part 7,
Max-Plank-Institut f\"ur Mathematik, Bonn,
 Preprint Series, 2007, 53, 1 -- 40.
 \bibitem[L.A.Gutnik,11(2007)L.A.Gutnik]{G11}
L.A.Gutnik, On some systems of difference equations, part 8,
Max-Plank-Institut f\"ur Mathematik, Bonn,
 Preprint Series, 2007, 64, 1 -- 44.
 \bibitem[L.A.Gutnik,12(2007)L.A.Gutnik]{G12}
L.A.Gutnik,, On some systems of difference equations, part 9,
Max-Plank-Institut f\"ur Mathematik, Bonn,
 Preprint Series, 2007, 129, 1 -- 36.
 \bibitem[L.A.Gutnik,13(2007)L.A.Gutnik]{G13}
 L.A.Gutnik, On some systems of
difference equations, part 10,
Max-Plank-Institut f\"ur Mathematik, Bonn, Preprint Series,2007,131, 1 -- 33.
 \bibitem[L.A.Gutnik,14(2007)L.A.Gutnik]{G14}
L.A.Gutnik, On some systems of difference equations, part 11,
Max-Plank-Institut f\"ur Mathematik, Bonn, Preprint Series, 2008,
 38, 1 -- 45.
 \bibitem[L.A.Gutnik,15(2006)L.A.Gutnik]{G15}
L.A.Gutnik,  On some systems of difference equations, Chebyshev Collection,
2006, v.7, No 3, , 140 -- 145.
\bibitem[L.A.Gutnik,16(2010)L.A.Gutnik]{G16}
L.A.Gutnik, Elementary Proof of Yu.V. Nesterenko expansion of the Nuber Zeta(3)
in Continued Fraction,Advances in Difference Equation, 2010,Articele Id
143521,11 pages.
\bibitem[L.A.Gutnik,17(2009)L.A.Gutnik]{G17}
L.A.Gutnik, On the number Zeta(3), Arxiv.org, Arxiv: 09022.4732.
 \bibitem[L.A.Gutnik,18(2013)L.A.Gutnik]{G18}
 On Expansion of Zeta(3) in Continued Fraction,
 (Short version),
International  Mathematical Forum, 2013, Vol.8, no 16, 771 -- 781.
   \end{thebibliography}
\end{document}